\documentclass[a4paper]{amsart}
\usepackage[utf8]{inputenc}
\usepackage{amsfonts}
\usepackage{amsmath,latexsym,amssymb,amsfonts, amsthm}
\usepackage{esint}
\usepackage{graphicx}
\usepackage{amscd}
\usepackage{color}
\usepackage{bm}             
\definecolor{my_red}{rgb}{0.8,0,0}
\usepackage{enumerate}
\usepackage{aliascnt}
\usepackage[dvips]{epsfig}
\usepackage{psfrag}
\usepackage{enumitem}
\usepackage{epsfig}
\usepackage{comment}

\usepackage[
  hmarginratio={1:1},     
  vmarginratio={1:1},     
  textwidth=16cm,        
  textheight=21cm,
  heightrounded,          
]{geometry}

\usepackage{graphicx,color}
\usepackage[colorlinks]{hyperref}
\usepackage{cleveref}
\hypersetup{linkcolor=blue,citecolor=blue,filecolor=black,urlcolor=blue}
\usepackage{mathtools}

\DeclarePairedDelimiter{\norm}{\lVert}{\rVert}

\renewcommand{\b}{\beta}

\definecolor{my_red}{RGB}{255, 0, 0}

\newtheorem{theorem}{Theorem}[section]
\crefname{theorem}{theorem}{theorems}
\Crefname{theorem}{Theorem}{Theorems}

\newaliascnt{proposition}{theorem}
\newtheorem{proposition}[proposition]{Proposition}
\aliascntresetthe{proposition}
\crefname{proposition}{proposition}{propositions}
\Crefname{proposition}{Proposition}{Propositions}

\newaliascnt{corollary}{theorem}

\aliascntresetthe{corollary}
\crefname{corollary}{corollary}{corollaries}
\Crefname{corollary}{Corollary}{Corollaries}

\newaliascnt{lemma}{theorem}
\newtheorem{lemma}[lemma]{Lemma}
\aliascntresetthe{lemma}
\crefname{lemma}{lemma}{lemmas}
\Crefname{lemma}{Lemma}{Lemmas}

\theoremstyle{definition}
\newaliascnt{remark}{theorem}
\newtheorem{remark}[remark]{Remark}
\aliascntresetthe{remark}
\theoremstyle{definition}
\crefname{remark}{remark}{remarks}
\Crefname{remark}{Remark}{Remarks}

\newaliascnt{example}{theorem}

\aliascntresetthe{example}
\theoremstyle{definition}
\crefname{example}{example}{examples}
\Crefname{example}{Example}{Examples}

\newaliascnt{definition}{theorem}

\aliascntresetthe{definition}
\theoremstyle{definition}
\crefname{definition}{definition}{definitions}
\Crefname{definition}{Definition}{Definitions}

\numberwithin{theorem}{section}

\numberwithin{equation}{section}

\newcommand{\beq}{\begin{equation}}
\newcommand{\eeq}{\end{equation}}

\newcommand{\R}{\mathbb{R}}

\newcommand{\mf}[1]{\mathbf{#1}}
\newcommand{\bs}[1]{\boldsymbol{#1}}

\renewcommand{\S}{\mathbb{S}}

\newcommand{\cA}{{\mathcal A}}   
\newcommand{\cB}{{\mathcal B}}   
\newcommand{\cC}{{\mathcal C}}   
   
\newcommand{\cE}{{\mathcal E}}

\newcommand{\cM}{{\mathcal M}}   
\newcommand{\cN}{{\mathcal N}}

\newcommand{\cU}{{\mathcal U}}

\usepackage{mathtools}

\DeclareMathOperator*{\argmax}{arg\,max}

\newcommand{\eps}{\varepsilon}

\DeclareMathOperator{\loc}{loc}

\usepackage{accents}
\newcommand{\ubar}[1]{\underaccent{\bar}{#1}}

\renewcommand{\epsilon}{\varepsilon}

\author[H. Fan, L. Giaretto, N. Soave]{Haining Fan, Lorenzo Giaretto and Nicola Soave}\thanks{}
\address{Haining Fan \newline \indent School of mathematics, China University of Mining and Technology,\newline \indent Xuzhou, Jiangsu 221116, China}
\email{fanhaining888@163.com}
\address{Lorenzo Giaretto and Nicola Soave\newline \indent
 Dipartimento di Matematica ``Giuseppe Peano'', Universit\`a di Torino, \newline \indent
Via Carlo Alberto 10,
10123 Torino, Italy}
\email{lorenzo.giaretto@unito.it, nicola.soave@unito.it}

\title[On a critical NLS system with $K$-wise interaction]{On a critical nonlinear Schr\"odinger system with $K$-wise interaction: the cooperative case}
\keywords{Nonlinear Schr\"odinger systems; $K$-wise interactions;
critical Sobolev exponent; cooperative systems;
least energy solutions}
\subjclass[2020]{35J50, 35B33; 35J47, 35B09}
\thanks{L.G. and N.S. are affiliated to the INDAM - GNAMPA group. L.G. was partially supported by the INDAM--GNAMPA 2026 research project ``Ottimizzazione Geometrica e Frontiere Libere''. H.F. is supported by the Fundamental Research Funds for the Central Universities (No. 2024KYJD2007)
and NSFC (No.12671144).\\
Declarations of interest: none. \\
Data availability: Data sharing not applicable to this article as no datasets were generated or analysed during the current study.}

\begin{document}

\maketitle

\begin{abstract}
We study least energy fully non-trivial solutions of a nonlinear Schr\"odinger system
with $K\ge3$ components and $K$-wise interaction, in the Sobolev critical regime. We consider bounded domains with
Brezis--Nirenberg-type assumptions, and the associated homogeneous problem
in $\mathbb R^N$.

For every positive coupling, we prove the existence of a nonnegative
ground state and identify a strictly positive threshold separating
semi-trivial ground states from fully non-trivial ones. When the
linear coefficients coincide, a finite-dimensional reduction gives
an exact formula for the ground state energy and a characterization
of all nonnegative ground states. In particular, above the threshold
all ground states are synchronized; for sufficiently strong
cooperation, their amplitude vector is unique and has an explicit
asymptotic behavior. We also establish existence results for positive
synchronized solutions beyond the ground state regime.

For sufficiently weak cooperation, we obtain least energy positive
solutions in the critical ternary case $N=K=3$, by minimizing on the
component-wise Nehari set and ruling out critical concentration.
The analysis includes the endpoint interaction exponent $q=1$,
where the energy functional is not differentiable.
\end{abstract}

\section{Introduction}

In this paper we study the existence and structure of least energy solutions
for a critical nonlinear Schr\"odinger system with simultaneous interaction
among all its components. More precisely, we consider
\begin{equation}\label{eq:system 1} 
\begin{cases}
-\Delta u_i+\lambda_i u_i
=\mu_i|u_i|^{2^*-2}u_i
+\beta f_q(u_i)\displaystyle\prod_{j\ne i}|u_j|^q
&\text{in }\Omega,\\
u_i\in H_0^1(\Omega), 
\end{cases} \qquad i=1,\ldots,K,
\end{equation}
where \(\Omega\) is either $\R^N$, or a smooth bounded domain,
\(\mu_i>0\), \(\beta \in \R\), and
\[
K\ge3,\qquad 2^*=\frac{2N}{N-2},\qquad q=\frac{2^*}{K}\ge1.
\]
Here \(f_q(s)=|s|^{q-2}s\) for \(q>1\), whereas
\(f_1(s)=\operatorname{sgn}(s)\), with \(\operatorname{sgn}(0)=0\).
The relation \(Kq=2^*\) makes both the self-interaction and the coupling
Sobolev critical. The assumptions \(K\ge3\) and \(q\ge1\) restrict the
admissible pairs to
\[
\begin{array}{c|c}
N&K\\ \hline
3&3,4,5,6\\
4&3,4\\
5&3\\
6&3
\end{array}
\]
and, in particular, \(1\le q\le2\). We focus on nonnegative solutions,
mainly in the cooperative regime $\beta>0$: the positive sign of \(\beta\) makes simultaneous
overlap of the components energetically favorable. The central questions are
when this attraction is sufficient to produce coexistence of all \(K\)
components, and what is the spatial structure of the resulting least energy solutions. The competitive case $\beta<0$ is also of interest, but requires completely different strategies and therefore is addressed in a separate companion paper \cite{FaGiSo26}, for readability.

Systems of coupled nonlinear Schr\"odinger equations are motivated,
among other contexts, by the study of interacting wave fields in
nonlinear optics. From this perspective, the components \(u_i\)
represent amplitudes of different fields, and the coefficients
\(\mu_i>0\) describe focusing self-interactions.
The product term models a collective interaction involving all
\(K\) components simultaneously, rather than a sum of pairwise
couplings. In the cooperative regime \(\beta>0\), its contribution in the energy functional favors the spatial overlap of all components, providing an energetic
mechanism for coexistence. The system can thus be viewed as an
idealized model for collective attractive interactions, and the
choice \(Kq=2^*\) places this mechanism in the Sobolev critical regime. It is worth
mentioning that systems with multiple (i.e., beyond binary) interaction terms naturally arise in various physical contexts; see, for
instance, \cite{CCO, LQZ, Pet, Pom} and references therein.

Let \(\lambda_1(\Omega)\) denote the first Dirichlet eigenvalue of
\(-\Delta\) in $\Omega$. Throughout the bounded-domain analysis, we assume that
\begin{equation}\label{coop-intro:linear-assumptions}
\begin{cases}
-\lambda_1(\Omega)<\lambda_i<0,&N\ge4,\\
-\lambda_1(\Omega)<\lambda_i<-\lambda^*(\Omega),&N=3,
\end{cases}
\qquad i=1,\ldots,K,
\end{equation}
where \(\lambda^*(\Omega)\) is the scalar Brezis--Nirenberg threshold
in dimension three. These hypotheses place each uncoupled equation in a
range where a positive scalar ground state exists and its energy lies
strictly below the corresponding whole-space Sobolev level
\cite{BN}. Thus \eqref{eq:system 1}  is a Brezis--Nirenberg-type
system in which the critical scalar equations are coupled through a single
\(K\)-fold product.

\medskip

Critical Schr\"odinger systems with binary interaction provide an important
point of comparison. The binary case has been studied extensively over the past few decades, starting with the seminal contributions \cite{AmCo, LiWe, MaMoPe, Si}. A comprehensive review of the literature is beyond the scope of this work. We therefore restrict our discussion to works directly relevant to the present study. In dimension four, Chen and Zou \cite{Chen_Zou_2012} studied
the system with cubic self-interactions and coupling terms
\(\beta uv^2\) and \(\beta vu^2\). They obtained least energy positive
solutions for sufficiently weak and sufficiently strong cooperation, and
classified the least energy solutions when the linear coefficients
coincide. In that symmetric linear setting, their results also exhibit
ranges of positive coupling where coexistence is impossible. Their subsequent
work \cite{Chen_Zou_2015} treats the critical binary system in dimensions
\(N\ge5\), where the exponent of each factor in the interaction is
\(2^*/2\in(1,2)\). In this range a least energy positive solution exists
for every \(\beta>0\), and synchronized least energy solutions can again
be constructed and classified for sufficiently strong cooperation. The
change from interaction exponent \(2\) to an exponent below \(2\) is
therefore already decisive in the binary problem.

The critical theory has also been developed for systems with an arbitrary
number of components coupled through sums of pairwise interactions,
\[
-\Delta u_i+\lambda_i u_i
=\mu_i|u_i|^{2^*-2}u_i
+\sum_{j\ne i}\beta_{ij}|u_i|^{r-2}u_i|u_j|^r,
\qquad r=\frac{2^*}{2}.
\]
Yin and Zou \cite{YZ22} established existence and classification results
for purely cooperative systems in higher dimensions. Tavares and You
\cite{TY20} analyzed critical systems in dimension four by organizing
components into groups and comparing the energy of the full system with
that of suitable subsystems. Tavares, You and Zou \cite{TYZ22} developed
this analysis in dimensions \(N\ge5\), including purely cooperative
systems and systems with mixed interactions. In particular, in the purely
cooperative case their ground states are positive for arbitrary positive
pairwise coupling coefficients. In dimension three, Liu, You and Zou
\cite{LYZ23} obtained least energy positive solutions under weak
cooperation. These works show how the exponent of the interaction, the
coupling coefficients and the energies of subsystems jointly govern
coexistence in critical pairwise models.

The coupling in \eqref{eq:system 1}  has a different structure. It
involves the product of all \(K\) components, rather than a sum of products
of two components. Two of the authors recently conducted a variational study of this simultaneous interaction in
the subcritical regime 
\cite{Giaretto_Soave_2026}. For systems with positive linear
coefficients and \(Kq<2^*\), we proved the existence of ground states
and identified a strictly positive threshold separating semi-trivial from
fully non-trivial ground states. We also obtained least energy positive
solutions for weak cooperation when \(q\ge2\). The present paper addresses the critical regime
\(Kq=2^*\), both on bounded domains and for the associated homogeneous
problem on \(\mathbb R^N\). Besides the loss of compactness at critical
growth, we also investigate questions concerning synchronized solutions and
the classification of critical ground states that were not addressed in
the subcritical analysis (although they make perfect sense in that case as well, see the later Remark \ref{rmk:synchronized subcritical case}).

Passing from binary to \(K\)-wise interaction changes the mechanism
responsible for coexistence; the usual infinitesimal perturbation excluding semi-trivial
ground states for arbitrarily weak cooperation is not available, unlike in the binary case for $N \ge 5$; we refer to the following Remark \ref{rem: k-wise vs binary} for a detailed discussion of this point.

Another structural difference concerns subsystems. If one component
vanishes identically, the entire $K$-wise interaction energy disappears and the
equations for the surviving components decouple into scalar
Brezis--Nirenberg problems. In a pairwise system, deleting a component leaves a coupled system among the remaining components. Thus the
inductive comparisons with lower-dimensional coupled subsystems used in
\cite{LYZ23, TY20,TYZ22, YZ22} cannot be transferred to the
\(K\)-wise problem. Here semi-trivial configurations are governed by sums of
scalar energies, whereas fully non-trivial ones are
governed by a genuinely collective interaction.

\subsection{Notation and main results}

We set \(H=H_0^1(\Omega,\mathbb R^K)\), denote by $|\cdot|_p$ the standard norm on $L^p(\R^d)$ or $L^p(\Omega)$, and
\begin{equation*}
\|v\|_{\lambda_i}^2
=\int_\Omega\bigl(|\nabla v|^2+\lambda_i v^2\bigr)dx.
\end{equation*}
The energy functional naturally associated with \eqref{eq:system 1} is the functional defined on $H$ as
\begin{equation}\label{energy}
E(\mathbf u)
=\frac12\sum_{i=1}^K\|u_i\|_{\lambda_i}^2
-\frac1{2^*}\sum_{i=1}^K\mu_i\int_\Omega|u_i|^{2^*}dx
-\frac\beta q \int_\Omega\prod_{i=1}^K|u_i|^qdx.
\end{equation}
A non-zero solution is called \emph{fully non-trivial} if
\(u_i\not\equiv0\) for every \(i\), and semi-trivial otherwise. A \emph{ground
state} is a minimizer of the energy among all non-zero solutions; a \emph{least energy
fully non-trivial solution} is a minimizer within the class of fully non-trivial solutions. These two notions need not coincide, even though the
system is cooperative.

Natural variational constraints to look for these solutions are the usual Nehari set and its component-wise counterpart,
\begin{align*}
\mathcal M'
&=\left\{\mathbf u\in H\setminus\{\mathbf0\}:
\sum_{i=1}^K\|u_i\|_{\lambda_i}^2
=\sum_{i=1}^K\mu_i\int_\Omega|u_i|^{2^*}dx
+K\beta \int_\Omega\prod_{i=1}^K|u_i|^qdx\right\},\\
\mathcal M
&=\left\{\mathbf u\in H:
u_i\not\equiv0,\quad
\|u_i\|_{\lambda_i}^2
=\mu_i\int_\Omega|u_i|^{2^*}dx+\beta \int_\Omega\prod_{i=1}^K|u_i|^qdx \ 
\text{  for every }i\right\}.
\end{align*}
Their associated minimal levels are
\begin{equation}\label{min levels}
\mathcal B'=\inf_{\mathcal M'}E,\qquad \mathcal B=\inf_{\mathcal M}E.
\end{equation}
Since \(Kq=2^*\), on either constraint 
\begin{equation*}
E(\mf{u}) = \frac1N \sum_{i=1}^K \|u_i\|_{\lambda_i}^2 \,;
\end{equation*}
moreover,
\(\mathcal M\subset\mathcal M'\), and hence \(\mathcal B'\le \mathcal B\). The level $\cB '$ is the ground state energy level, whereas $\cB$ is the
least energy level among fully non-trivial configurations satisfying
the component-wise Nehari constraints. The first minimization allows semi-trivial configurations, while the
second imposes \(K\) separate identities, and, if $\mathcal{B}$ is achieved, it is necessary to show that the constraint is natural in the
regime under consideration.

\medskip

Note that, for $q>1$, the functional $E$ is of class $C^1$ on $H$, and its
nonnegative critical points are precisely the weak solutions of
\eqref{eq:system 1}. Moreover, the identities defining $\cM'$ and $\cM$ are equivalent to
\[
E'(\mf{u})[\mf{u}]=0 \quad \text{and} \quad E'(\mf u)[0,\dots,0,u_i,0,\dots,0]=0,
\qquad i=1,\dots,K,
\]
respectively. In particular, any ground state belongs to $\cM'$, any fully non-trivial solution of \eqref{eq:system 1}
belongs to $\mathcal M$, and it is classical to show that $\mathcal{M}'$ is a natural constraint: any constrained critical point of $E$ on $\mathcal{M}'$ is in fact a critical point of $E$ in the full space $H$.

A word of caution is needed when $q=1$, since the interaction term is not differentiable on the whole space $H$. Whenever no
differentiability is required, we will nevertheless use \eqref{energy}
as the energy associated with the original nonnegative problem and
define the corresponding Nehari sets directly through their explicit identities. On the other hand, when differentiability is needed, we introduce a suitable auxiliary functional and corresponding constraints to overcome this lack of regularity (see the definitions of $\widetilde E$ and $\widetilde{\mathcal{M}}'$ in \eqref{eq:E tilde def} and \eqref{eq:M tilde prime def} below).

We now state the main results of the paper. We first consider the case of identical linear coefficients: in addition to \eqref{coop-intro:linear-assumptions}, we suppose that
\begin{equation}\label{coop-intro:lambda symmetric assumption}
	\lambda_1=\cdots=\lambda_K=\lambda.
\end{equation}
In this case, it is natural to investigate the existence of positive \emph{synchronized solutions}, namely solutions of
the form \((\ell_1\omega,\ldots,\ell_K\omega)\), where \(\omega\)
is a positive scalar solution of
\begin{equation}\label{BN}
-\Delta\omega+\lambda\omega=\omega^{2^*-1},
\qquad \omega\in H_0^1(\Omega),
\end{equation}
and the positive amplitudes solve
\begin{equation}\label{eq:algebraic-system}
\mu_i\ell_i^{2^*-2}
+\beta\ell_i^{q-2}\prod_{j\ne i}\ell_j^q=1,
\qquad i=1,\ldots,K.
\end{equation}
Note that, under assumption \eqref{coop-intro:linear-assumptions}, a positive solution $\omega$ of \eqref{BN} always exists, and in particular one can take a positive ground state of \eqref{BN}, see \cite{BN}.

Let
$
\mathbb R_{>0}^K
:=
\left\{
(x_1,\dots,x_K)\in\mathbb R^K:
x_i>0\text{ for every }i=1,\dots,K
\right\}.
$
Without loss of generality, we assume that
\[
\mu_1\geq\mu_2\geq\cdots\geq\mu_K>0.
\]
In our first result, we provide sufficient conditions for existence of synchronized solutions. 
\begin{theorem}\label{thm:existence synchronized solutions}
	Assume \eqref{coop-intro:linear-assumptions}, \eqref{coop-intro:lambda symmetric assumption} and let $\omega$ be a least energy positive solution of \eqref{BN}. Then the following statements hold.
	\begin{enumerate}
		\item If $1<q<2$, then there exists
		$(\ell_1,\dots,\ell_K)\in\R_{>0}^K$ such that $(\ell_1\omega,\dots,\ell_K\omega)$
		is a positive synchronized solution of \eqref{eq:system 1}.
		
		\item If $q=2$ (this happens when $N=3$, $K=3$) and
		\[
		0<\beta<
		(\mu_2\mu_3)^{\frac{1}{2}}
		\qquad\text{or}\qquad
		\beta>
		(\mu_1\mu_3)^{\frac{1}{2}},
		\]
		then there exists
		$(\ell_1,\ell_2,\ell_3)\in\R_{>0}^3$ such that $(\ell_1\omega,\ell_2 \omega,\ell_3\omega)$ is a positive synchronized solution of \eqref{eq:system 1}.
		\item If $q=1$ and
		\[
		0<\beta<
		(\mu_2\cdots\mu_{K-1})^{\frac{1}{K-2}}
		\qquad\text{or}\qquad
		\beta>
		(\mu_1\mu_4\cdots\mu_K)^{\frac{1}{K-2}},
		\]
		then there exists
		$(\ell_1,\dots,\ell_K)\in\R_{>0}^K$ such that $(\ell_1\omega,\dots,\ell_K\omega)$	is a positive synchronized solution of \eqref{eq:system 1}.
	\end{enumerate}
\end{theorem}

\begin{remark}\label{rem:synchronized_ranges}
	\begin{enumerate}
		\item In the proof of case $q=2$ (see Section \ref{sec:synchronized}), if $\mu_1=\mu_2$, one can actually obtain existence for all $\beta>0$, thereby filling the gap $\beta=(\mu_1\mu_3)^{\frac{1}{2}}$. See Remark \ref{rmk:mu1 = mu2} below.
		\item In the parameter ranges of $\beta$ excluded in items (2) and (3) of Theorem~\ref{thm:existence synchronized solutions}, our approach does not establish nonexistence; whether positive synchronized solutions exist in these ranges remains open.
		\item In fact, for $q=1$, depending on the values of $\mu_i$, the upper range for $\beta$ may overlap with the lower range, yielding the existence of positive synchronized solutions for every $\beta>0$.
	\end{enumerate}
\end{remark}




\begin{remark}
Let $\omega$ be a positive  solution of
\begin{equation*}
-\Delta\omega=\omega^{2^*-1},
\qquad \omega \in D^{1,2}(\mathbb{R}^N).
\end{equation*}
Then a result similar to Theorem~\ref{thm:existence synchronized solutions}  can be 
obtained for \eqref{eq:system 1} defined on the domain $\mathbb{R}^N$ instead of $\Omega$. 
\end{remark}

We next investigate the variational properties of synchronized solutions. 
It turns out that there is a sharp threshold above which all ground states are synchronized. 

\begin{theorem}\label{thm:synchronized are least energy for beta large}
	Assume \eqref{coop-intro:lambda symmetric assumption}. Then there exists $\beta_0>0$ (defined in \eqref{eq:beta_0 def})
	such that, for every $\beta >\beta_0$, all nonnegative ground states are positive and synchronized. 
\end{theorem}

Regarding the limit case $\beta=\beta_0$, we refer to Remark \ref{ref: on beta_0 611} for more details.

Moreover, we describe the asymptotic behavior of the amplitude vectors $\bs \ell(\b)$ and, for sufficiently strong cooperation, establish a uniqueness result.

\begin{theorem}\label{thm:uniqueness and asymptotics}
	Assume \eqref{coop-intro:lambda symmetric assumption}. Then
    $$
    \b^{1/(2^*-2)}\ell_i(\b)\to 1 \quad \text{as } \b\to+\infty.
    $$
    Moreover, there exists
	$\beta_1>0$ such that, if
	$\beta\geq\beta_1$, there is a unique $\bs \ell(\b)$ such that $(\ell_1(\b)\omega, \dots, \ell_K(\b)\omega)$ is a ground state solution. As a consequence, the synchronized ground state is unique if $\Omega\subset\mathbb{R}^N$ is a ball and $\b\geq\b_1$.
\end{theorem}

\begin{remark}\label{rmk:synchronized subcritical case}
The results of Theorems \ref{thm:existence synchronized solutions}, \ref{thm:synchronized are least energy for beta large}, and \ref{thm:uniqueness and asymptotics} can be readily extended to the subcritical setting considered in \cite{Giaretto_Soave_2026}, provided that a positive solution to the corresponding scalar problem exists. In particular, this condition is satisfied for
\[
-\Delta \omega+\lambda\omega=\omega^{p-1}
\quad\text{in }\mathbb{R}^N,
\]
when \(\lambda>0\) and \(2<p<2^*\), and for
\[
\begin{cases}
-\Delta \omega+\lambda\omega=\omega^{p-1} & \text{in }\Omega,\\
\omega=0 & \text{on }\partial\Omega,
\end{cases}
\]
when \(\Omega\subset\mathbb{R}^N\) is a bounded domain, \(\lambda>-\lambda_1(\Omega)\), and \(2<p<2^*\).
\end{remark}

In the proof of the previous results, we observe that the algebraic system for the amplitude vector reduces to scalar equations linked by a common interaction parameter. Each equation has two positive branches on a suitable parameter interval, and the branch choices must satisfy a compatibility condition. The analysis of this condition yields the existence ranges stated above and is carried out in Section \ref{sec:synchronized}.

We now turn to the general case of possibly distinct linear coefficients. Under assumption \eqref{coop-intro:linear-assumptions}, we aim to investigate the existence of ground states and of least energy fully non-trivial
solutions.

\begin{theorem}\label{thm:ground-state-bounded-cooperative}
	Let $\Omega\subset\mathbb R^N$ be a bounded domain, assume that
	$N\geq 3$, and let \eqref{coop-intro:linear-assumptions} hold.
	Then, for every $\beta>0$, the level $\cB '$ defined in \eqref{min levels} is attained by a
	nonnegative ground state of \eqref{eq:system 1}. Moreover, there exists $\bar\beta>0$ (explicitly defined in \eqref{eq:def bar beta bounded} below) such that:
	\begin{itemize}
		\item if $\beta>\bar\beta$, every nonnegative ground state is
		positive. In particular,
		\[
		\cB=\cB ',
		\]
		and $\cB$ is attained;
		
		\item if $0<\beta<\bar\beta$, no ground state is fully non-trivial.
	\end{itemize}
\end{theorem}

\begin{remark}
Assume that \(\lambda_1=\cdots=\lambda_K\) and \(1<q<2\).
By Theorem~\ref{thm:existence synchronized solutions}-(1), positive synchronized solutions exist for every
\(\beta>0\). However, Theorem~\ref{thm:ground-state-bounded-cooperative} shows that, for
\(0<\beta<\bar\beta\), every ground state is semi-trivial.
In particular, the positive synchronized solutions in this range
have energy strictly greater than \(\mathcal{B}'\).
Thus, the threshold \(\bar\beta\) describes a transition in the
structure of ground states, and need not be a threshold for the
existence of positive solutions. This distinction is essential
when studying coexistence through energy minimization.
\end{remark}
\begin{remark}
In the symmetric case $\lambda_1=\dots=\lambda_K$, the threshold $\bar\beta$ coincides with the value $\beta_0$ introduced in Theorem \ref{thm:synchronized are least energy for beta large}, as required for consistency between the two results. We refer to Remark \ref{rmk:beta_0 equal beta bar} below for the proof.
\end{remark}

For sufficiently weak cooperation, minimizing \(E\) on
\(\mathcal M'\) yields only semi-trivial solutions.
To obtain coexistence in this regime, we turn to minimization
on the component-wise Nehari set \(\mathcal M\).
The next result establishes attainment of \(\cB\) in the critical
ternary case \(N=K=3\), where \(q=2\).
Compared with the subcritical result in \cite{Giaretto_Soave_2026}, this requires
overcoming the loss of compactness at critical growth. 

\begin{theorem}\label{thm:weak cooperative bounded domain}
    Let $\Omega\subset\mathbb R^3$ be a bounded regular domain, assume that
	$K= 3$, and let \eqref{coop-intro:linear-assumptions} hold. Then, there exists $\ubar \b>0$ such that, if $0<\b<\ubar\beta$, $\cB$ is attained by a least energy positive solution of \eqref{eq:system 1}.
\end{theorem}

\begin{remark}\label{rem: k-wise vs binary}
Passing from binary to \(K\)-wise interaction changes the mechanism
responsible for coexistence,
as illustrated by the following heuristic
argument. Let \(w\) be a positive ground state of the first scalar
equation, and perturb the semi-trivial configuration
\(\mathbf u_0=(w,0,\ldots,0)\) for system \eqref{eq:system 1} by setting
\[
\mathbf u_\varepsilon
=(w,\varepsilon\varphi_2,\ldots,\varepsilon\varphi_K),
\]
where the nonnegative profiles \(\varphi_j\) are fixed and $w \cdot \varphi_2 \cdots \varphi_K \not \equiv 0$. Introducing the missing components produces a positive quadratic
contribution to the energy, of order \(\varepsilon^2\), whereas their
scalar nonlinearities contribute only at order \(\varepsilon^{2^*}\).
The cooperative interaction lowers the energy by a term of order
\(\beta\varepsilon^{(K-1)q}\), since all \(K-1\) missing components
must be activated simultaneously for the product to become nonzero.
This balance is preserved by the common rescaling
\(t_\varepsilon\mathbf u_\varepsilon\) onto the global Nehari set:
\(t_\varepsilon\to1\) as $\eps \to 0$, and, for fixed \(\beta>0\),
\[
E(t_\varepsilon\mathbf u_\varepsilon)-E(\mathbf u_0)
=C_1\varepsilon^2
-\beta C_2\varepsilon^{(K-1)q}
+o\bigl(\varepsilon^2+\varepsilon^{(K-1)q}\bigr),
\]
where
\[
C_1=\frac12\sum_{j=2}^K\|\varphi_j\|_{\lambda_j}^2>0,
\qquad
C_2=\frac1q\int_\Omega w^q\prod_{j=2}^K\varphi_j^q \; >0.
\]
We note that \((K-1)q\ge2\) for every choice of $K$ and $N$ and, if the inequality is strict, the quadratic cost dominates for every fixed \(\beta\) as
\(\varepsilon\to0\); in the borderline case, the sign depends on the
coefficients of the two quadratic contributions.
For a binary interaction with exponent \(r\in(1,2)\), the corresponding
gain has instead order \(\varepsilon^r\), and therefore dominates the quadratic
cost for every positive coupling. 
Therefore, the usual infinitesimal perturbation excluding semi-trivial
ground states for arbitrarily weak cooperation is no longer available.
This local comparison motivates the strictly positive coexistence
threshold $\bar \beta$ established in Theorem \ref{thm:ground-state-bounded-cooperative}, in contrast with the fact that this threshold vanishes for the binary case considered in \cite{Chen_Zou_2015}. 
\end{remark}


Since both the nonlinearities and the coupling term in \eqref{eq:system 1}
have critical growth, the existence of least energy solutions is closely
related to the corresponding limiting problem in the whole space,
\begin{equation}\label{eq:system D12}
	\begin{cases}
		-\Delta u_i
		=
		\mu_i |u_i|^{2^\ast-2}u_i
		+
		\beta f_q(u_i)
		\displaystyle\prod_{\substack{j=1\\j\neq i}}^K |u_j|^q
		& \text{in }\mathbb R^N,\\[0.4em]
		u_i\in D^{1,2}(\mathbb R^N),
	\end{cases} \qquad i=1,\dots,K,
\end{equation}
where, as before, $\mu_i>0$, $q=2^\ast/K$, and
\[
D^{1,2}(\mathbb R^N)
:=
\left\{
u\in L^{2^\ast}(\mathbb R^N):
|\nabla u|\in L^2(\mathbb R^N)
\right\}, \quad \|u\|_{D^{1,2}}
:=
\left(
\int_{\mathbb R^N}|\nabla u|^2\,dx
\right)^{1/2}.
\]
We denote by $S$ the sharp Sobolev constant for the embedding
$D^{1,2}(\mathbb R^N)\hookrightarrow L^{2^\ast}(\mathbb R^N)$, namely
\begin{equation}\label{eq:sharp Sobolev R^N}
	\int_{\mathbb R^N}|\nabla u|^2\,dx
	\geq
	S
	\left(
	\int_{\mathbb R^N}|u|^{2^\ast}\,dx
	\right)^{\frac{2}{2^\ast}}
	\qquad
	\text{for every }u\in D^{1,2}(\mathbb R^N).
\end{equation}
Moreover, for $\varepsilon>0$ and $y\in\mathbb R^N$, we denote by
\[
U_{\varepsilon,y}(x)
:=
[N(N-2)]^{\frac{N-2}{4}}
\left(
\frac{\varepsilon}
{\varepsilon^2+|x-y|^2}
\right)^{\frac{N-2}{2}}
\]
the Aubin--Talenti bubble. The family $\left\{
U_{\varepsilon,y}:
\varepsilon>0,\ y\in\mathbb R^N
\right\}$ exhausts all positive solutions of
\[
-\Delta u=|u|^{2^\ast-2}u
\qquad\text{in }\mathbb R^N,
\]
see \cite{CaGiSp}; each $U_{\varepsilon,y}$ attains equality in \eqref{eq:sharp Sobolev R^N}, and in particular
\begin{equation*}
	\int_{\mathbb R^N}|\nabla U_{\varepsilon,y}|^2\,dx
	=
	\int_{\mathbb R^N}|U_{\varepsilon,y}|^{2^\ast}\,dx
	=
	S^{N/2}.
\end{equation*}

System \eqref{eq:system D12} always admits semi-trivial
solutions, for instance those of type
\begin{equation}\label{eq:semitrivial solution D12}
    \mf U_i
:=
\left(
0,\dots,0,
\mu_i^{-\frac{N-2}{4}}U_{\varepsilon,y},
0,\dots,0
\right),
\end{equation}
and also those with more than one
non-zero component. Our main interest is therefore in fully non-trivial
solutions and in determining when ground states are fully non-trivial. Set $\mathcal D
:=
D^{1,2}(\mathbb R^N,\mathbb R^K)$, and define
\begin{equation}\label{eq:energy-D12}
	I(\mf u)
	=
	\frac12
	\sum_{i=1}^K
	\int_{\mathbb R^N}|\nabla u_i|^2\,dx
	-
	\frac{1}{2^\ast}
	\int_{\mathbb R^N}
	\left(
	\sum_{i=1}^K\mu_i|u_i|^{2^\ast}
	+
	K\beta\prod_{i=1}^K|u_i|^q
	\right)\,dx.
\end{equation}
For $q>1$, $I\in C^1(\mathcal D)$. When $q=1$, we shall use
\eqref{eq:energy-D12} as the energy of the original nonnegative
problem, without claiming differentiability on the whole space.

As in the bounded-domain setting, we introduce the usual Nehari set $\mathcal{N}'$ and its component-wise counterpart $\mathcal{N}$, together with the associated minimal levels $\mathcal{A}':=\inf_{\mathcal{N}'} I$ and $\mathcal{A}:=\inf_{\mathcal{N}} I$.

The following result describes the behavior of the energy levels and will play a key role in the study of the problem on bounded domains.

\begin{theorem}\label{thm:ground-state-D12}
	Let $N\geq3$, $K\geq3$, and $q=2^\ast/K\geq1$. Then,
	\begin{enumerate}[label=(\roman*)]	
		\item for every $\b>0$, $\cA'$ is attained by a nonnegative
		ground state of \eqref{eq:system D12}.
		\item there exists $\beta_0>0$ (defined in \eqref{eq:beta_0 def}) such that
		\begin{itemize}
			\item if $\beta>\beta_0$, every nonnegative ground state is
			positive and is of the form
            \begin{equation}\label{eq:ground state synchronized form}
                (a_1 U_{\varepsilon,y},\dots, a_K U_{\varepsilon,y})
            \end{equation}
            for some $a_i>0$, $\varepsilon>0$ and $y\in \R^N$.
            Hence $\cA=\cA'$;
			\item if $0<\beta<\beta_0$, every ground state is
			semi-trivial and is of the form \eqref{eq:semitrivial solution D12} with $i\in \argmax_j \mu_j$ and some $\eps>0$ and $y\in \R^N$.
		\end{itemize}
	\end{enumerate}
\end{theorem}

\begin{remark}
    We emphasize that the threshold $\beta_0$ in the theorem above coincides with that of Theorem~\ref{thm:synchronized are least energy for beta large}. Indeed, the underlying algebraic problem determining the nature of ground states in both settings is identical, allowing the proofs to be presented in a unified manner. We refer to Section~\ref{sec:auxiliary} for details.
\end{remark}

\begin{remark}
    The results concerning the existence of synchronized solutions in Theorem~\ref{thm:existence synchronized solutions} extend seamlessly to the problem on the entire space $\mathbb{R}^N$ by replacing the positive solutions on bounded domains with Aubin--Talenti bubbles. In particular, synchronized solutions continue to exist for $0 < \beta < \beta_0$ (at least when $1 < q < 2$), though by the preceding theorem they do not constitute ground states. However, in the specific case $q=2$ (namely $N = K = 3$), we prove below that for sufficiently small $\beta$ they do indeed form least energy positive solutions.
\end{remark}

\begin{theorem}\label{thm:weak cooperative D12}
If $q=2$ (namely $N=3$ and $K=3$) and $\b>0$ is small enough, then $\cA$ is attained by a least energy positive solution of the form \eqref{eq:ground state synchronized form}.
\end{theorem}

The whole-space analysis is more than an auxiliary existence result:
it identifies the energy cost of concentration. For the ground-state
problem on a bounded domain, we prove the strict inequality \(\mathcal{B}'<\mathcal{A}'\),
which prevents a minimizing sequence from disappearing entirely into a
critical bubble. 

In the weakly cooperative triple-interaction problem,
minimization on \(\mathcal M\) requires further comparisons. Let
\(b_i\) be the scalar ground-state energy of the \(i\)-th equation on
\(\Omega\), and set
\(c_i=N^{-1}\mu_i^{-(N-2)/2}S^{N/2}\). We establish
\[
\cB<\cA,\qquad \cB<b_i+b_j+c_k
\quad\text{for every }\{i,j,k\}=\{1,2,3\}.
\]
These estimates exclude both total concentration and partial loss of
components. The latter comparisons reflect the \(K\)-wise structure:
once a component is absent, the surviving equations are scalar, and a
missing component can carry the energy of a scalar Sobolev bubble.
The argument also uses a Brezis--Lieb-type splitting
\cite{BL83} adapted to the full product, so that the critical interaction
is accounted for together with the separate scalar nonlinearities. The case $q<2$ seems even harder and remains open.

\subsection*{Structure of the paper}
We conclude the introduction with a brief description of the organization
of the paper. Section~\ref{sec:synchronized} is devoted to synchronized
solutions: we study the associated algebraic system and prove Theorem \ref{thm:existence synchronized solutions}. In Section~\ref{sec:auxiliary},
we treat in a unified way the existence and characterization of ground states in the symmetric case on bounded domains and on the whole $\R^N$, proving Theorems \ref{thm:synchronized are least energy for beta large}, \ref{thm:uniqueness and asymptotics} and \ref{thm:ground-state-D12}. These results
provide the main variational ingredients for the bounded-domain problem.
In Section~\ref{sec:ground states bounded domain}, we consider the cooperative case on
bounded domains and prove Theorem~\ref{thm:ground-state-bounded-cooperative}.
Finally, Section~\ref{sec:weak cooperative} is devoted to the weak cooperative case
and contains the proofs of Theorems~\ref{thm:weak cooperative bounded domain} and \ref{thm:weak cooperative D12}.

\section{Existence of synchronized solutions}\label{sec:synchronized}

\subsection{Proof of Theorem \ref{thm:existence synchronized solutions}}

The first part of the section is devoted to the proof of Theorem \ref{thm:existence synchronized solutions}. We first show the existence results.
\begin{proof}[Proof of Theorem \ref{thm:existence synchronized solutions}.] Recall that finding synchronized solutions reduces to finding solutions
of the algebraic system \eqref{eq:algebraic-system}. Set
	\[
	T:=\beta\prod_{j=1}^K\ell_j^q.
	\]
	Then \eqref{eq:algebraic-system} is equivalent to
	\begin{equation*}
		\ell_i^2-\mu_i\ell_i^{2^*}=T,
		\qquad i=1,\dots,K.
	\end{equation*}
	For every $i=1,\dots,K$, define
	\[
	f_i(\ell):=\ell^2-\mu_i\ell^{2^*}
	=
	\ell^2\left(1-\mu_i\ell^{2^*-2}\right),
	\qquad \ell>0.
	\]
	Setting	$b_i:=\mu_i^{-\frac{1}{2^*-2}}$,	we have
	\[
	f_i(\ell)>0
	\quad\Longleftrightarrow\quad
	0<\ell<b_i.
	\]
	Moreover, $f_i$ attains its unique maximum at
	\[
	\widehat\ell_i
	:=
	\left(\frac{2}{2^*\mu_i}\right)^{\frac{1}{2^*-2}},
	\]
	and
	\[
	T_i
	:=
	f_i(\widehat\ell_i)
	=
	\frac{2^*-2}{2^*}
	\left(\frac{2}{2^*\mu_i}\right)^{\frac{2}{2^*-2}}.
	\]
	Since $\mu_1\geq\cdots\geq\mu_K$, we have
	\[
	\overline T
	:=
	\min_{1\leq i\leq K}T_i
	=
	T_1.
	\]
	For every $T\in(0,\overline T)$, there exist unique
	$\ell_i^-(T),\ell_i^+(T)>0$ such that
	\[
	f_i(\ell_i^-(T))
	=
	f_i(\ell_i^+(T))
	=
	T,
	\qquad
	0<\ell_i^-(T)<\widehat\ell_i<\ell_i^+(T)<b_i.
	\]
	Furthermore,
	\begin{equation}\label{min levelsranches-limit}
		\ell_i^-(T)\to0,
		\qquad
		\ell_i^+(T)\to b_i
		\qquad\text{as }T\to0^+,
	\end{equation}
	and, since
	\[
	T=(\ell_i^-(T))^2\left(1-\mu_i(\ell_i^-(T))^{2^*-2}\right),
	\]
	we have
	\begin{equation}\label{eq:lminus-asymptotics}
		\ell_i^-(T)\sim T^{1/2}
		\qquad\text{as }T\to0^+.
	\end{equation}
	Therefore, it remains to choose one of the two branches for each
	component in such a way that
	\begin{equation*}
		\beta=T \Big(\prod_{i=1}^K\ell_i(T)^q\Big)^{-1}.
	\end{equation*}
	
	Now, assume at first that $1<q<2$. Define
	\[
	\Phi_1(T)
	:=
	\begin{cases}
		\bigl(
		\ell_1^+(T),\ell_2^-(T),
		\ell_3^+(T),\dots,\ell_K^+(T)
		\bigr),
		& T\in[0,\overline T],\\[0.3em]
		\bigl(
		\ell_1^-(2\overline T-T),
		\ell_2^-(2\overline T-T),
		\ell_3^+(2\overline T-T),\dots,
		\ell_K^+(2\overline T-T)
		\bigr),
		& T\in[\overline T,2\overline T],
	\end{cases}
	\]
	and
	\[
	g_1(T)
	:=
	\begin{cases}
		\displaystyle
		\frac{T}
		{\prod_{i=1}^K\Phi_{1,i}(T)^q},
		& T\in[0,\overline T],\\[1.8em]
		\displaystyle
		\frac{2\overline T-T}
		{\prod_{i=1}^K\Phi_{1,i}(T)^q},
		& T\in[\overline T,2\overline T].
	\end{cases}
	\]
    By construction, $g_1\in C([0,2\bar T))$, and by \eqref{min levelsranches-limit} and \eqref{eq:lminus-asymptotics},
	\[
	g_1(T)\to0
	\quad\text{as }T\to0^+,
	\qquad
	g_1(T)\to+\infty
	\quad\text{as }T\to2\overline T^-,
	\]
	where we have used $1<q<2$. Hence, by continuity, for every $\beta>0$
	there exists $\widehat T\in(0,2\overline T)$ such that
	$g_1(\widehat T)=\beta$. Thus $\Phi_1(\widehat T)$ solves
	\eqref{eq:algebraic-system}.
	
	Assume now that $q=2$. Consider first the path
	\[
	\Phi_2(T)
	:=
	\begin{cases}
		\bigl(
		\ell_1^+(T),\ell_2^+(T),\ell_3^+(T)
		\bigr),
		& T\in[0,\overline T],\\[0.3em]
		\bigl(
		\ell_1^-(2\overline T-T),
		\ell_2^+(2\overline T-T),
		\ell_3^+(2\overline T-T)
		\bigr),
		& T\in[\overline T,2\overline T],
	\end{cases}
	\]
	and define $g_2$ accordingly as before. Then $g_2(T)\to0$ as $T\to0^+$, while
	\[
	g_2(T)
	\to
	(b_2b_3)^{-2}
	=
	(\mu_2\mu_3)^{\frac{1}{2}}
	\quad\text{as }T\to2\overline T^-.
	\]
	Therefore a solution exists whenever
	\[
	0<\beta<
	(\mu_2\mu_3)^{\frac{1}{2}}.
	\]
	On the other hand, considering the path
	\[
	\Phi_3(T)
	:=
	\begin{cases}
		\bigl(
		\ell_1^+(T),\ell_2^-(T),
		\ell_3^+(T)
		\bigr),
		& T\in[0,\overline T],\\[0.3em]
		\bigl(
		\ell_1^-(2\overline T-T),
		\ell_2^-(2\overline T-T),
		\ell_3^+(2\overline T-T)
		\bigr),
		& T\in[\overline T,2\overline T],
	\end{cases}
	\]
	the corresponding function $g_3$ satisfies
	\[
	g_3(T)
	\to
	(b_1b_3)^{-2}
	=
	(\mu_1\mu_3)^{\frac{1}{2}}
	\quad\text{as }T\to0^+,
	\]
	and
	\[
	g_3(T)\to+\infty
	\quad\text{as }T\to2\overline T^-.
	\]
	Hence a solution exists also whenever
	\[
	\beta>
	(\mu_1\mu_3)^{\frac{1}{2}}.
	\]
	
	Finally, assume that $q=1$. Consider
	\[
	\Phi_4(T)
	:=
	\begin{cases}
		\bigl(
		\ell_1^-(T),\ell_2^-(T),\ell_3^-(T),
		\ell_4^+(T),\dots,\ell_K^+(T)
		\bigr),
		& T\in[0,\overline T],\\[0.3em]
		\bigl(
		\ell_1^+(2\overline T-T),
		\ell_2^-(2\overline T-T),
		\ell_3^-(2\overline T-T),
		\ell_4^+(2\overline T-T),\dots,
		\ell_K^+(2\overline T-T)
		\bigr),
		& T\in[\overline T,2\overline T].
	\end{cases}
	\]
	The corresponding function $g_4$ satisfies $g_4(T)\to+\infty$ as $T\to0^+$, and
	\[
	g_4(T)
	\to
	(b_1b_4\cdots b_K)^{-1}
	=
	(\mu_1\mu_4\cdots\mu_K)^{\frac{1}{K-2}}
	\quad\text{as }T\to2\overline T^-.
	\]
	Thus a solution exists whenever
	\[
	\beta>
	(\mu_1\mu_4\cdots\mu_K)^{\frac{1}{K-2}}.
	\]
	Similarly, along the path
	\[
	\Phi_5(T)
	:=
	\begin{cases}
		\bigl(
		\ell_1^+(T),\dots,\ell_{K-1}^+(T),\ell_K^-(T)
		\bigr),
		& T\in[0,\overline T],\\[0.3em]
		\bigl(
		\ell_1^-(2\overline T-T),
		\ell_2^+(2\overline T-T),\dots,
		\ell_{K-1}^+(2\overline T-T),
		\ell_K^-(2\overline T-T)
		\bigr),
		& T\in[\overline T,2\overline T],
	\end{cases}
	\]
	the corresponding function $g_5$ satisfies $g_5(T)\to0$ as $T\to0^+$, and
	\[
	g_5(T)
	\to
	(b_2\cdots b_{K-1})^{-1}
	=
	(\mu_2\cdots\mu_{K-1})^{\frac{1}{K-2}}
	\quad\text{as }T\to2\overline T^-.
	\]
    And so a solution exists whenever
    \[
	0<\beta<
	(\mu_2\cdots\mu_{K-1})^{\frac{1}{K-2}}.
	\]
	This concludes the proof.
\end{proof}

\begin{remark}\label{rmk:mu1 = mu2}
    In the case $q=2$, if $\mu_1=\mu_2$ the present proof would address the existence of a synchronized solution for every $\b\neq(\mu_1 \mu_2)^{1/2}$. In fact, in this case, we can define
    \[
	\Phi(T)
	:=
	\begin{cases}
		\bigl(
		\ell_1^+(T),\ell_2^+(T),\ell_3^+(T)
		\bigr),
		& T\in[0,\overline T],\\[0.3em]
		\bigl(
		\ell_1^-(2\overline T-T),
		\ell_2^-(2\overline T-T),
		\ell_3^+(2\overline T-T)
		\bigr),
		& T\in[\overline T,2\overline T],
	\end{cases}
	\]
    so that the corresponding $g(t)$ satisfies $g(t)\to 0$ as $t\to 0^+$ and $g(t)\to +\infty$ as $t\to 2\bar T$, proving the existence of a synchronized solution for every $\b>0$.
\end{remark}

\section{Synchronized ground states in the symmetric setting}\label{sec:auxiliary}
The present section is devoted to the proof of the results concerning the existence and classification of ground states, both for the problem on bounded domains under the symmetric assumption \eqref{coop-intro:lambda symmetric assumption} and on $\R^N$. The proofs of Theorems \ref{thm:synchronized are least energy for beta large}, \ref{thm:uniqueness and asymptotics}, and \ref{thm:ground-state-D12} are unified within an abstract setting that applies to both cases simultaneously.

We consider either of the following two settings:
\begin{equation}\label{note:settings}
\begin{array}{c|c|c|c|c}
 D & X & Q_X(v) & L_X & \text{Problem} \\ \hline
 \Omega & H_0^1(\Omega) & \displaystyle\int_\Omega (|\nabla v|^2+\lambda v^2)dx & -\Delta+\lambda & \eqref{eq:system 1} \\[5pt]
 \R^N & D^{1,2}(\R^N) & \displaystyle\int_{\R^N}|\nabla v|^2dx & -\Delta & \eqref{eq:system D12}
\end{array}
\end{equation}
where, in the bounded-domain case, \(\Omega\) is smooth, bounded and connected, and assumptions \eqref{coop-intro:linear-assumptions} and \eqref{coop-intro:lambda symmetric assumption} hold. In both cases, the quadratic form \(Q_X\) is positive definite, and we set
\begin{equation}\label{note:scalarconstant}
 S_X = \inf_{v\in X\setminus\{0\}} \frac{Q_X(v)}{\|v\|_{L^{2^*}(D)}^2},
\end{equation}
which is positive and attained.

A positive scalar ground state \(\omega\) is normalized so that
\[
 L_X\omega = \omega^{2^*-1}, \qquad Q_X(\omega) = \int_D\omega^{2^*}dx = S_X^{N/2}.
\]
In the $D^{1,2}(\R^N)$ case, the positive scalar ground states are the normalized
Aubin--Talenti functions \(U_{\varepsilon,y}\).

For \(\mf u\in X^K\), set
\begin{align*}
 Q(\mf u)&=\sum_{i=1}^K Q_X(u_i),\\
 P_\beta(\mf u)&=\int_D\left(\sum_{i=1}^K\mu_i|u_i|^{2^*}
                  +K\beta\prod_{i=1}^K|u_i|^q\right)dx,\\
 E_X(\mf u)&=\frac12Q(\mf u)-\frac{1}{2^*}P_\beta(\mf u),\\
 \cN'_X&=\{\mf u\in X^K\setminus\{0\}:Q(\mf u)=P_\beta(\mf u)\},
 \qquad \cC'_{X}=\inf_{\mf u\in\cN'_X}E_X(\mf u).
\end{align*}
These definitions apply also when \(q=1\), without differentiating
\(E_X\).

\medskip

For every $\bs a\in \mathbb S^{K-1}$ with $a_i\geq 0$ for every $i$, define
\begin{equation}\label{eq:M_beta def}
 F_\beta(\bs a)=\sum_{i=1}^K\mu_i a_i^{2^*}
                  +K\beta\prod_{i=1}^K a_i^q,
 \qquad
 M_\beta=\max_{\substack{a_i\ge0\\|\bs a|=1}}F_\beta(\bs a).
\end{equation}

We first relate this algebraic problem to the ground state energy.

\begin{proposition}\label{prop:connection between ground state energy and algebraic problem}
In either setting \eqref{note:settings}, for every \(\beta>0\),
\begin{equation*}
 \cC'_{X}=\frac{1}{N} S_X^{N/2} M_\beta^{-\frac{N-2}{2}}.
\end{equation*}
This level is attained by a nonnegative solution of
the corresponding problem, and hence is the ground state energy.
Moreover, all nonnegative ground states, and only those, have the form
\begin{equation*}
 \mf u=M_\beta^{-\frac{N-2}{4}}\bs a\,\omega,
 \qquad \bs a\in\argmax_{\substack{a_i\ge0\\|\bs a|=1}}F_\beta,
\end{equation*}
where \(\omega\) is a positive scalar ground state.
\end{proposition}

\begin{proof}
\noindent \textbf{Step 1:} $\cC_X'\geq \frac{1}{N}S_X^{N/2} M_\beta^{-(N-2)/2}$.

First, since $Kq = 2^*$, $P_\beta$ is homogeneous of degree $2^*$. Thus, whenever $\mf u \neq \bs 0$ and $P_\beta(\mf u) > 0$, there exists a unique scalar $t_{\mf u} > 0$ such that $t_{\mf u}\mf u \in \cN'_X$, given explicitly by
\[
t_{\mf u} = \left( \frac{Q(\mf u)}{P_\beta(\mf u)} \right)^{\frac{1}{2^*-2}} = \left( \frac{Q(\mf u)}{P_\beta(\mf u)} \right)^{\frac{N-2}{4}}.
\]
Consequently, the ground state energy level $\cC'_X = \inf_{\cN'_X} E_X$ can be characterized as
\begin{equation}\label{eq:abstract-quotient}
\cC'_{X} = \frac{1}{N} \left( \inf_{\substack{\mf u \in X^K \setminus \{\bs 0\} \\ P_\beta(\mf u) > 0}} \frac{Q(\mf u)}{P_\beta(\mf u)^{2/2^*}} \right)^{N/2}.
\end{equation}

Since $\beta > 0$, clearly $M_\beta>0$. For any $\mf u = (u_1,\dots,u_K) \in X^K \setminus \{\bs 0\}$, let 
\[
w := \left( \sum_{i=1}^K u_i^2 \right)^{1/2} \in X.
\]
At every point where $w(x) > 0$, applying the definition of $M_\beta$ to the unit vector 
\[
\frac{1}{w(x)}(|u_1(x)|,\dots,|u_K(x)|) \in \mathbb S^{K-1}
\]
yields
\[
\sum_{i=1}^K \mu_i |u_i|^{2^*} + K\beta \prod_{i=1}^K |u_i|^q \le M_\beta w^{2^*} \qquad \text{a.e. in } D.
\]
Integrating over $D$, we obtain
\begin{equation}\label{eq:P-bound-abstract}
P_\beta(\mf u) \le M_\beta \int_D w^{2^*}dx.
\end{equation}
Moreover,
\begin{equation}\label{eq:4}
|\nabla w|^2\leq \sum_{i=1}^K|\nabla u_i|^2\qquad\text{a.e. in } D.
\end{equation} 
In the bounded-domain case, since $\lambda \sum_i u_i^2 = \lambda w^2$, the linear term satisfies $Q(\mf u) \ge Q_X(w)$ independently of the sign of $\lambda$. Combining this with \eqref{note:scalarconstant}, we find
\[
Q(\mf u) \ge Q_X(w) \ge S_X \left( \int_D w^{2^*} dx \right)^{2/2^*} \ge S_X M_\beta^{-2/2^*} P_\beta(\mf u)^{2/2^*}.
\]
Taking the infimum over all $\mf u$ with $P_\beta(\mf u) > 0$ and recalling \eqref{eq:abstract-quotient}, we find
\begin{equation}\label{eq:sharp-bound-abstract}
\cC'_{X} \ge \frac{1}{N} S_X^{N/2} M_\beta^{-\frac{N-2}{2}}.
\end{equation}

\medskip

\noindent\textbf{Step 2:} The above inequality is an equality and it is attained by a ground state.

Let $\bs a = (a_1,\dots,a_K)\in \mathbb S^{K-1}$ be a nonnegative maximizer for $M_\beta$, and let $\omega \in X$ be a positive scalar ground state normalized such that 
\[
Q_X(\omega) = \int_D \omega^{2^*}dx = S_X^{N/2}.
\]
Setting $\mf U := (a_1\omega,\dots,a_K\omega)$, we have
\[
Q(\mf U) = Q_X(\omega) = S_X^{N/2}, \qquad P_\beta(\mf U) = M_\beta \int_D \omega^{2^*}dx = M_\beta S_X^{N/2}.
\]
Therefore, $\mf U$ attains equality in \eqref{eq:sharp-bound-abstract}. Multiplying $\mf U$ by $t_{\mf U} = M_\beta^{-(N-2)/4}$ yields 
\[
\mf V := t_{\mf U}\mf U = M_\beta^{-(N-2)/4}\bs a\,\omega \in \cN'_X,
\]
which achieves the infimum, thus establishing $\cC'_X = \frac{1}{N} S_X^{N/2} M_\beta^{-(N-2)/2}$.

If $q > 1$, $\cN'_X$ is a natural constraint for the $C^1$ functional $E_X$, so $\mf V$ is a nonnegative critical point and hence a ground state solution of the corresponding problem in \eqref{note:settings}.

If $q = 1$, the nonlinearity in $E_X$ contains non-differentiable absolute values. To overcome this, we introduce the smooth auxiliary functional $\widetilde{E}_X: X^K \to \mathbb R$ defined by
\[
\widetilde{E}_X(\mf u) := \frac{1}{2} Q(\mf u) - \frac{1}{2^*} \int_D \left( \sum_{i=1}^K \mu_i |u_i|^{2^*} + 2^* \beta \prod_{i=1}^K u_i \right) dx.
\]
Defining the corresponding auxiliary Nehari manifold
\[
\widetilde{\cN}'_X := \left\{ \mf u \in X^K \setminus \{\bs 0\} : Q(\mf u) = \int_D \left( \sum_{i=1}^K \mu_i |u_i|^{2^*} + 2^* \beta \prod_{i=1}^K u_i \right) dx \right\},
\]
we observe that for every $\mf u \in X^K$ and $\beta > 0$, $\widetilde{E}_X(|\mf u|) \le \widetilde{E}_X(\mf u)$, where $|\mf u| := (|u_1|,\dots,|u_K|)$. Consequently,
\begin{equation*}
 \cC'_X = \inf_{\cN'_X} E_X = \inf_{\substack{\mf u \in \widetilde{\cN}'_X \\ u_i \ge 0}} \widetilde{E}_X(\mf u) = \inf_{\widetilde{\cN}'_X} \widetilde{E}_X.
\end{equation*}
Since $\mf V = M_\beta^{-(N-2)/4}\bs a\,\omega$ has nonnegative components $V_i \ge 0$, it satisfies $E_X(\mf V) = \widetilde{E}_X(\mf V)$ and $\mf V \in \widetilde{\cN}'_X$. Because $\widetilde{\cN}'_X$ is a natural constraint for the $C^1$ functional $\widetilde{E}_X$, $\mf V$ is a weak solution of
\begin{equation*}
L_X V_i = \mu_i |V_i|^{2^*-2}V_i + \beta \prod_{j \neq i} V_j \qquad \text{in } D, \quad i = 1,\dots,K.
\end{equation*}
Standard elliptic regularity implies that $\mf V$ is in $W^{2,p}_{\loc}$ for some $p\geq 1$, thus $L_X V_i=0$ a.e. on $\{V_i = 0\}$, so
\[
\prod_{j \neq i} V_j = 0 \qquad \text{a.e. on } \{V_i = 0\}.
\]
Thus, $\mf V$ is a ground state solution to the corresponding problem (\eqref{eq:system 1} and \eqref{eq:system D12}, respectively).

\noindent\textbf{Step 3:} characterization of ground states in terms of maximizers of $M_\b$.

Let $\mf u$ be an arbitrary nonnegative ground state. Then all inequalities in the chain leading to \eqref{eq:sharp-bound-abstract} must hold as equalities. In particular,
\[
Q_X(w) = S_X \left( \int_D w^{2^*}dx \right)^{2/2^*},
\]
so $w$ is a positive scalar ground state profile in $D$. Set $a_i := u_i / w$, then
$|\bs a| = 1$. The pointwise identity
\[
\sum_{i=1}^K |\nabla u_i|^2 - |\nabla w|^2 = w^2 \sum_{i=1}^K \left| \nabla \left( \frac{u_i}{w} \right) \right|^2
\]
and equality in \eqref{eq:4} imply that $\nabla(u_i/w) = 0$ a.e. in $D$. Since $D$ is connected, the ratios $a_i$ are constant functions. Thus $\mf u = \bs a w$ for some constant vector $\bs a \in \mathbb S^{K-1}$ with $a_i \ge 0$.
Finally, equality in \eqref{eq:P-bound-abstract} forces $F_\beta(\bs a) = M_\beta$, proving that $\bs a$ is a nonnegative maximizer of $M_\beta$, and this concludes the proof.
\end{proof}

Set $\S^{K-1}_+:=\{\bs a\in \S^{K-1}: a_i>0 \text{ for every }i=1,\dots,K\}$. We now introduce the threshold that characterizes the nature of maximizers of \eqref{eq:M_beta def}. 

\begin{proposition}\label{prop:beta_0}
Define
\begin{equation}\label{eq:beta_0 def}
 \beta_0:=\frac1K \inf_{\bs a\in \S^{K-1}_+}
 \frac{\mu_{\max}-\displaystyle\sum_{i=1}^K\mu_i a_i^{2^*}}{\displaystyle\prod_{i=1}^K a_i^q},
\end{equation}
where $\mu_{\max}:=\displaystyle\max_{1\leq i\leq K}\mu_i$.
Then $\beta_0>0$. Moreover,
\begin{itemize}
    \item if \(0<\beta<\beta_0\), the maximizers in \eqref{eq:M_beta def} are exactly the vertices
    \(e_j\) with \(\mu_j=\mu_{\max}\);
    \item if \(\beta>\beta_0\), all maximizers have all strictly positive coordinates.
\end{itemize}
\end{proposition}

\begin{proof}
We first show that $\beta_0>0$. Let
\[
G(\bs a) := \frac{\mu_{\max}-\displaystyle\sum_{i=1}^K\mu_i a_i^{2^*}}{\displaystyle\prod_{i=1}^K a_i^q}.
\]
First, notice that for every $\bs a\in \mathbb S^{K-1}_+$ we have
\[
\sum_{i=1}^K\mu_i a_i^{2^*} \leq \mu_{\max}\sum_{i=1}^K a_i^{2^*} < \mu_{\max} \Big(\sum_{i=1}^K a_i^2\Big)^{2^*/2} = \mu_{\max},
\]
thus $G(\bs a)>0$ for every $\bs a\in\mathbb S^{K-1}_+$.

It remains to exclude that $\beta_0=0$. Suppose by contradiction that there exists $\bs a_n\in\mathbb S^{K-1}_+$ such that $G(\bs a_n)\to0$. Up to a subsequence, $\bs a_n\to\bs a$ for some $\bs a$ such that $|\bs a|=1$ and $a_i\geq 0$ for every $i=1,\dots,K$. If $\bs a\in\mathbb S^{K-1}_+$, by continuity we would have $G(\bs a_n)\to G(\bs a)>0$, a contradiction. Thus $a_i=0$ for at least one index $i$, and consequently
\[
\prod_{i=1}^K a_{i,n}^q\to0.
\]
Since $G(\bs a_n)\to0$, necessarily
\[
\mu_{\max} - \sum_{i=1}^K\mu_i a_{i,n}^{2^*} \to0 \quad \implies \quad 
\sum_{i=1}^K\mu_i a_i^{2^*}=\mu_{\max}.
\]
On the other hand,
\[
\sum_{i=1}^K\mu_i a_i^{2^*} \leq \mu_{\max}\sum_{i=1}^K a_i^{2^*} \leq \mu_{\max} \left(\sum_{i=1}^K a_i^2\right)^{2^*/2} = \mu_{\max},
\]
and equality is possible only if $\bs a$ has exactly one nonzero component, say $a_j$, such that $\mu_j=\mu_{\max}$. Set
\[
\delta_n^2 :=\sum_{i\neq j}a_{i,n}^2\to0.
\]
Then $a_{j,n}=(1-\delta_n^2)^{1/2}$, and therefore
\begin{equation*}
\mu_{\max} -\sum_{i=1}^K\mu_i a_{i,n}^{2^*} = \mu_{\max} \left[ 1-(1-\delta_n^2)^{2^*/2} \right] - \sum_{i\neq j}\mu_i a_{i,n}^{2^*} = \frac{2^*}{2}\mu_{\max}\delta_n^2 +o(\delta_n^2),
\end{equation*}
where in the last equality we used that
\[
\sum_{i\neq j}a_{i,n}^{2^*} \leq \Big(\sum_{i\neq j}a_{i,n}^2\Big)^{2^*/2} = \delta_n^{2^*} = o(\delta_n^2).
\]
In particular, for $n$ sufficiently large,
\[
\mu_{\max} -\sum_{i=1}^K\mu_i a_{i,n}^{2^*} \geq c\delta_n^2
\]
for some $c>0$. Moreover,
\[
\prod_{i=1}^K a_{i,n}^q = a_{j,n}^q\prod_{i\neq j}a_{i,n}^q \leq C\delta_n^{q(K-1)},
\]
since $a_{i,n}\leq\delta_n$ for every $i\neq j$. Since $q\geq1$ and $K\geq3$, $q(K-1)\geq2$. Consequently,
\[
G(\bs a_n) \geq C\delta_n^{2-q(K-1)}.
\]
If $q(K-1)=2$, the right-hand side is bounded from below by a positive constant, whereas if $q(K-1)>2$, it tends to $+\infty$. In either case this contradicts $G(\bs a_n)\to0$. Therefore $\beta_0>0$.

Suppose now that $\beta<\beta_0$. Then, for every $\bs a\in\mathbb S^{K-1}$ with $a_i>0$ for all $i$, the definition of $\beta_0$ gives
\[
\sum_{i=1}^K\mu_ia_i^{2^\ast} + K\beta\prod_{i=1}^Ka_i^q < \mu_{\max}.
\]
On the other hand, if $i_0$ is such that $\mu_{i_0}=\mu_{\max}$, then the coordinate vector $\bs e_{i_0}$ satisfies
\[
\sum_{i=1}^K\mu_i(\bs e_{i_0})_i^{2^\ast} + K\beta\prod_{i=1}^K(\bs e_{i_0})_i^q = \mu_{\max}.
\]
This concludes the proof in this case.

Suppose instead that $\beta>\beta_0$. By the definition of $\beta_0$, there exists $\bs a\in\mathbb S^{K-1}_+$ such that
\[
\sum_{i=1}^K\mu_ia_i^{2^\ast} + K\beta\prod_{i=1}^Ka_i^q > \mu_{\max}.
\]
Consequently, $M_\beta>\mu_{\max}$. On the other hand, if $\bs b\in\mathbb S^{K-1}$ has at least one zero component, then the interaction term vanishes and
\[
\sum_{i=1}^K\mu_i|b_i|^{2^\ast} \leq \mu_{\max}\sum_{i=1}^K|b_i|^{2^\ast} \leq \mu_{\max} \left(\sum_{i=1}^K b_i^2\right)^{2^*/2} = \mu_{\max}.
\]
Therefore no maximizer of $M_\beta$ can have a zero component and every nonnegative maximizer of $M_\beta$ lies in $\mathbb S^{K-1}_+$. This completes the proof.
\end{proof}

Combining the last two propositions, we prove Theorem \ref{thm:synchronized are least energy for beta large} and Theorem \ref{thm:ground-state-D12}.
\begin{proof}[Proof of Theorems \ref{thm:synchronized are least energy for beta large} and \ref{thm:ground-state-D12}]
    The results follow immediately combining Proposition \ref{prop:connection between ground state energy and algebraic problem} and \ref{prop:beta_0}.
\end{proof}

\begin{remark}\label{ref: on beta_0 611}
    Notice that the case $\beta=\beta_0$ is genuinely excluded from Theorem \ref{thm:synchronized are least energy for beta large}: whenever the infimum in
\eqref{eq:beta_0 def} is attained at some $\boldsymbol a_*\in\mathbb{S}^{K-1}_+$ (which is
always the case when $q(K-1)>2$, since then $G(\boldsymbol a)\to+\infty$ as
$\boldsymbol a\to\partial\mathbb{S}^{K-1}_+$) one has that
$F_{\beta_0}(\boldsymbol a_*)=\mu_{\max}=F_{\beta_0}(\boldsymbol e_{i_0})$ for every
$i_0$ with $\mu_{i_0}=\mu_{\max}$, so that for $\beta=\beta_0$ both positive
synchronized and semi-trivial ground states coexist.
\end{remark}

\begin{lemma}\label{lem:unique algebraic maximizer}
For all sufficiently large \(\beta\), \(F_\beta\) has a unique
nonnegative maximizer \(\bs a_\beta\). The corresponding amplitude
vector \(\bs \ell_\beta=M_\beta^{-(N-2)/4}\bs a_\beta\) satisfies
\begin{equation}\label{note:asymptotic}
 \beta^{(N-2)/4}\ell_i(\beta)\to 1,
 \qquad i=1,\ldots,K.
\end{equation}
\end{lemma}

\begin{proof}[Proof]
We begin by proving \eqref{note:asymptotic}. Consider the normalized functional given by
\[
 \frac{F_\beta(\bs a)}{K\beta} = \prod_{i=1}^K a_i^q + \frac{1}{K\beta} \sum_{i=1}^K \mu_i a_i^{2^*}.
\]
Clearly its maximizers among nonnegative elements in $\S^{K-1}$ are the same of $M_\b$. Since $\sum_{i=1}^K \mu_i a_i^{2^*}$ is continuous and thus bounded on the compact set $\mathbb S^{K-1}_+$, we have
\[
 \frac{F_\beta(\bs a)}{K\beta} \longrightarrow \prod_{i=1}^K a_i^q \qquad \text{uniformly on } \mathbb S^{K-1}_+ \text{ as } \beta \to +\infty.
\]
By the arithmetic-geometric mean inequality, for any $\bs a \in \mathbb S^{K-1}_+$,
\[
 \prod_{i=1}^K a_i^q = \left( \prod_{i=1}^K a_i^2 \right)^{q/2} \le \left( \frac{1}{K} \sum_{i=1}^K a_i^2 \right)^{Kq/2} = \left( \frac{1}{K} \right)^{2^*/2} = K^{-2^*/2},
\]
and equality holds if and only if $a_1^2 = \dots = a_K^2 = 1/K$. This shows that
\[
 \bs a_* = \left( \frac{1}{\sqrt{K}}, \dots, \frac{1}{\sqrt{K}} \right)
\]
is the unique maximizer of the limit problem
$$
\max_{\S^{K-1}_+} \prod_{i=1}^K a_i^q.
$$
By uniform convergence on the compact space $\mathbb S^{K-1}_+$, any sequence of maximizers $\bs a_\beta \in \argmax_{\mathbb S^{K-1}_+} F_\beta$ must converge to $\bs a_*$ as $\beta \to +\infty$. In particular, the maximum value satisfies
\[
 \frac{M_\beta}{K\beta} = \max_{\mathbb S^{K-1}_+} \frac{F_\beta(\bs a)}{K\beta} \to \prod_{i=1}^K (a_*)_i^q = K^{-2^*/2} \qquad \text{as } \beta \to +\infty.
\]
Recall that the ground state amplitudes are given by $\bs \ell_\beta = M_\beta^{-(N-2)/4} \bs a_\beta$, the last equation yields
\[
 \beta^{\frac{N-2}{4}} \ell_i(\beta) = \beta^{\frac{N-2}{4}} M_\beta^{-\frac{N-2}{4}} a_{i,\beta} \sim \beta^{\frac{N-2}{4}} \left( \beta K^{1-2^*/2} \right)^{-\frac{1}{2^*-2}} a_{i,\beta} = K^{\frac{2^*-2}{2(2^*-2)}} a_{i,\beta} = K^{1/2} a_{i,\beta}.
\]
Since $a_{i,\beta} \to (a_*)_i = K^{-1/2}$, \eqref{note:asymptotic} is proved.

\medskip

We now turn to the proof of uniqueness. The amplitudes $\ell_i$ satisfy the algebraic system
\[
 \ell_i^2 = \mu_i \ell_i^{2^*} + \beta \prod_{j=1}^K \ell_j^q, \qquad i = 1, \dots, K.
\]
Expressed in terms of $r_i = \beta^{(N-2)/4} \ell_i$ and $\varepsilon = \beta^{-1} \to 0^+$, using the identity $\frac{2^*(N-2)}{4} = \frac{N}{2}$, the above equation becomes
\[
 \beta^{-\frac{N-2}{2}} r_i^2 = \varepsilon \mu_i \beta^{-\frac{N-2}{2}} r_i^{2^*} + \beta^{1 - \frac{N}{2}} \prod_{j=1}^K r_j^q.
\]
Noting that $1 - \frac{N}{2} = -\frac{N-2}{2}$, multiplying the equation by $\beta^{\frac{N-2}{2}}$ cancels all powers of $\beta$ except in the perturbation term, yielding:
\[
 G_i(\boldsymbol r, \varepsilon) := r_i^2 - \varepsilon \mu_i r_i^{2^*} - \prod_{j=1}^K r_j^q = 0, \qquad i = 1, \dots, K.
\]
The map $\bs G = (G_1, \dots, G_K) : \mathbb R^K \times \mathbb R \to \mathbb R^K$ is smooth in a neighborhood of $(\boldsymbol r, \varepsilon) = (\bs 1, 0)$. Clearly $G_i(\bs 1,0)=0$ for every $i$. We now evaluate the Jacobian matrix $D_{\boldsymbol r} \bs G(\mathbf{1}, 0)$. Computing the partial derivatives gives
\[
 \frac{\partial \bs G_i}{\partial r_k}(\mathbf{1}, 0) = 2 \delta_{ik} - q \prod_{j \neq k} 1^{q} \cdot 1^{q-1} = 2 \delta_{ik} - q.
\]
In matrix form,
\[
 D_{\boldsymbol r} G(\mathbf{1}, 0) = 2 \operatorname{Id} - q J,
\]
where $J = (J_{ij})_{i,j=1}^K$ is the rank-one matrix with $J_{ij} = 1$ for every $i,j = 1, \dots, K$.
We now compute the spectrum of $D_{\boldsymbol r} G(\mathbf{1}, 0)$:
\begin{itemize}
 \item For any vector $\bs v \in \mathbb R^K$ such that $\sum_{i=1}^K v_i = 0$ (a space of dimension $K-1$), we have $J \bs v = \bs 0$, so
 \[
  D_{\boldsymbol r} G(\mathbf{1}, 0) \bs v = 2 \bs v.
 \]
 Thus, $\lambda_1 = 2$ is an eigenvalue with multiplicity $K-1$.
 \item For the vector $\mathbf{1} = (1, \dots, 1)^{\mathsf T}$, we have $J \mathbf{1} = K \mathbf{1}$, so
 \[
  D_{\boldsymbol r} G(\mathbf{1}, 0) \mathbf{1} = (2 - Kq) \mathbf{1} = (2 - 2^*) \mathbf{1}.
 \]
 Thus, $\lambda_2 = 2 - 2^* < 0$ is a simple eigenvalue (since $2^* > 2$).
\end{itemize}
Since all eigenvalues are non-zero, $D_{\boldsymbol r} G(\mathbf{1}, 0)$ is strictly invertible. By the Implicit Function Theorem, there exist $\varepsilon_0 > 0$ and a neighborhood $U \subset \mathbb R^K$ of $\mathbf{1}$ such that for every $\varepsilon \in [0, \varepsilon_0)$, the equation $G(\boldsymbol r, \varepsilon) = \bs 0$ has a \emph{unique} solution $\boldsymbol r(\varepsilon)$ in $U$.

Assume now by contradiction that there are two sequences of maximizers $\bs a_\b\neq \bs b_\b$ for $M_\b$ as $\b\to\infty$: by \eqref{note:asymptotic}, the rescaled vectors 
$$
\boldsymbol r_\beta=\b^{(N-2)/4} M_\b^{-(N-2)/4} \bs a_\b, \quad \boldsymbol r'_\beta=\b^{(N-2)/4} M_\b^{-(N-2)/4} \bs b_\b
$$ 
both satisfy $\boldsymbol r_\beta\to \bs 1, \bs r'_\b \to \bs 1$ as $\beta \to +\infty$. Thus, for $\beta$ sufficiently large, violates the uniqueness of solution to $G(\bs r,1/\beta)=0$ in $U$. The proof is therefore complete.
\end{proof}

\begin{remark}
Lemma~\ref{lem:unique algebraic maximizer} gives uniqueness of the amplitude vector of
ground states. It does not assert uniqueness of all positive solutions
of the algebraic system. On \(\Omega\), the remaining freedom is the
choice of the scalar ground state; if that scalar state is unique,
so is the positive least energy solution of the system. In the
whole space, the remaining parameters are the common translation
and dilation of the Aubin--Talenti profile.
\end{remark}

\begin{proof}[Proof of Theorem \ref{thm:uniqueness and asymptotics}]
    It follows from Lemma \ref{lem:unique algebraic maximizer}. If the domain is a ball, since the positive ground state solution of the scalar Brezis-Nirenberg problem is unique (see \cite{AdimurthiYadava1994}), the ground state is unique.
\end{proof}

\section{Ground states on bounded domains}\label{sec:ground states bounded domain}

This section contains the proof of Theorem \ref{thm:ground-state-bounded-cooperative}. We set
    \begin{equation*}
        \bar S:=\min_{1\leq i\leq K}\inf_{u\in H_0^1(\Omega)\setminus\{0\}}\frac{\norm{u}_{\lambda_i}^2}{\mu_i^{2/2^*}\norm{u}_{2^*}^2},
    \end{equation*}
the optimal constant associated with the Brezis--Nirenberg quotients in $H_0^1(\Omega)$ for the coefficients $\lambda_i$. As pointed out in the introduction, by
\cite{BN}, the Brezis--Nirenberg problem
\[
-\Delta u+\lambda_i u
=
\mu_i|u|^{2^*-2}u,
\qquad
u\in H_0^1(\Omega),
\]
has a positive least energy solution
$u_{\mu_i}\in C^2(\Omega)\cap C(\overline{\Omega})$ with energy
\begin{equation}\label{eq:B_mu def least energy bounded domain}
	B_{\mu_i}:=
	\frac12\int_\Omega
	\left(
	|\nabla u_{\mu_i}|^2
	+\lambda_i u_{\mu_i}^2
	\right)dx
	-
	\frac{\mu_i}{2^*}
	\int_\Omega |u_{\mu_i}|^{2^*}dx
	<\frac1N
	\mu_i^{\frac{2-N}{2}}S^{N/2},
	\qquad i=1,\dots,K.
\end{equation} 

We first introduce a crucial constant that characterizes the nature of
ground states. The proof follows exactly the same argument as
\cite[Theorem 1.1 and Proposition 1.3]{Giaretto_Soave_2026}.

\begin{proposition}\label{prop:beta bar charact H01}
	Let $N\geq3$, $K\geq3$, $q=2^*/K\geq1$, and define
	\begin{equation}\label{eq:def bar beta bounded}
	   \bar\beta
	:=
	\inf
	\left\{
	\frac{
		\bar S^{-\frac{2^*}{2}}
		\left(
		\sum_{i=1}^K\norm{u_i}_{\lambda_i}^2
		\right)^{\frac{2^*}{2}}
		-
		\sum_{i=1}^K\mu_i\norm{u_i}_{2^*}^{2^*}
	}{
		K\displaystyle\int_\Omega
		\prod_{i=1}^K|u_i|^q\,dx
	}
	:
	\mf u\in H_0^1(\Omega,\R^K),\
	\prod_{i=1}^K u_i\not\equiv0
	\right\},    
	\end{equation}
    Then $\bar\beta>0$. Moreover,
	\begin{enumerate}
		\item if $0<\beta<\bar\beta$, then any ground state is
		semi-trivial, with only one non-trivial component;
		
		\item if $\beta>\bar\beta$, then any ground state is fully
		non-trivial and, up to a change of sign of some components,
		is a least energy positive solution of \eqref{eq:system 1}.
	\end{enumerate}
\end{proposition}

\begin{remark}\label{rmk:beta_0 equal beta bar}
    In the symmetric case $\lambda_1=\dots=\lambda_K=\lambda$, we can show that $\bar\b=\beta_0$, where $\beta_0$ is introduced in \eqref{eq:beta_0 def}. Indeed, first observe that in that case
    $$
    \bar S=\mu_{\max}^{-2/2^*} \frac{\norm{\omega}_{\lambda}^2}{\norm{\omega}_{2^*}^2},
    $$
    where $\omega$ is the least energy positive solution of $-\Delta \omega+\lambda \omega= \omega^{2^*-1}$. Then, denoting by 
\[
S_\lambda:= \inf_{v \neq 0} \frac{\|v\|_\lambda^2}{\|v\|_{2^*}^2},
\]
for every $\mf u\in H_0^1(\Omega,\R^K)$ with $\prod_i u_i \neq 0$ we have 
    \begin{align*}
    &\frac{
		\mu_{\max} S_\lambda^{-\frac{2^*}{2}}
		\left(
		\sum_{i=1}^K\norm{u_i}_{\lambda}^2
		\right)^{\frac{2^*}{2}}
		-
		\sum_{i=1}^K\mu_i\norm{u_i}_{2^*}^{2^*}
	}{
		K\displaystyle\int_\Omega
		\prod_{i=1}^K|u_i|^q\,dx
	}\geq \frac{
		\mu_{\max} S_\lambda^{-\frac{2^*}{2}}
		\left(
		S_\lambda\sum_{i=1}^K\norm{u_i}_{2^*}^2
		\right)^{\frac{2^*}{2}}
		-
		\sum_{i=1}^K\mu_i\norm{u_i}_{2^*}^{2^*}
	}{
		K \prod_{i=1}^K\norm{u_i}_{2^*}^q
	} \\
    & \qquad \geq \inf_{\bs a\in \R^K_{>0}} \frac{
		\mu_{\max} \left(\sum_{i=1}^K a_i^2\right)^{\frac{2^*}{2}}
		-
		\sum_{i=1}^K\mu_i a_i^{2^*}
	}{
		K \prod_{i=1}^K a_i^q
	} = \inf_{\bs a\in \S^{K-1}_{+}} \frac{
		\mu_{\max} -\sum_{i=1}^K\mu_i a_i^{2^*}
	}{
		K \prod_{i=1}^K a_i^q
	}=\b_0,
    \end{align*}
    thus $\bar\beta \geq \beta_0$. On the other hand, choosing as test function $u_i=a_i \omega$, with $a_i>0$ for every $i$, yields
    \begin{equation*}
        \bar \beta \leq\frac{
		\mu_{\max} S_\lambda^{-\frac{2^*}{2}}
		\left(
		\sum_{i=1}^Ka_i^2 \norm{\omega}_{\lambda}^2
		\right)^{\frac{2^*}{2}}
		-
		\sum_{i=1}^K\mu_i a_i^{2^*}\norm{\omega}_{2^*}^{2^*}
	}{
		K\left(\prod_{i=1}^K a_i^q\right)\norm{\omega}_{2^*}^{2^*}
	}= \frac{
		\mu_{\max} \left(\sum_{i=1}^K a_i^2\right)^{\frac{2^*}{2}}
		-
		\sum_{i=1}^K\mu_i a_i^{2^*}
	}{
		K \prod_{i=1}^K a_i^q
	}
    \end{equation*}
    Taking the infimum over $\bs a\in \R^K_{>0}$, we obtain the reverse inequality, thus $\bar\beta=\beta_0$.
\end{remark}

First, we need to establish some preliminary lemmas.
\begin{lemma}
	Assume that \eqref{coop-intro:linear-assumptions} holds, and let $N\geq 3$. Then the following holds:
	\begin{equation}\label{min levelsprime-less-Aprime}
		\cB '<\cA'.
	\end{equation}
\end{lemma}
\begin{proof}
	By Theorem \ref{thm:ground-state-D12}, there exists
	$\bs a\in\R^K\setminus\{\bs0\}$, with $a_i\geq0$ for every
	$i=1,\dots,K$, such that $\mf U=(a_1U_{1,0},\dots,a_KU_{1,0})\in\cN'$ is a ground state of problem \eqref{eq:system D12}.	Set
	\[
	\lambda_{\bs a}
	:=
	\frac{\sum_{i=1}^K\lambda_i a_i^2}
	{\sum_{i=1}^K a_i^2}.
	\]
	Since $\lambda_{\bs a}$ is a convex combination of the $\lambda_i$'s,
	our assumptions imply
    \[
	-\lambda_1(\Omega)<\lambda_{\bs a}<0
	\quad \text{if }N\geq4, \quad \text{and} \quad -\lambda_1(\Omega)
	<
	\lambda_{\bs a}
	<
	-\lambda^*(\Omega)
	\quad\text{if }N=3.
    \]
	Hence, by the Brezis--Nirenberg result (see \cite[Lemmas 1.1, 1.3]{BN} for the case $N=3$ and $N\geq 4$, respectively),
	\begin{equation*}
		S_{\lambda_{\bs a}}(\Omega)
		:=
		\inf_{w\in H_0^1(\Omega)\setminus\{0\}}
		\frac{\displaystyle\int_\Omega
			\left(
			|\nabla w|^2+\lambda_{\bs a}w^2
			\right)\,dx
		}{		\displaystyle	\left(
			\int_\Omega|w|^{2^*}\,dx
			\right)^{2/2^*}
		}
		<S,
	\end{equation*}
	where $S$ denotes the optimal Sobolev constant in $\R^N$.
	Therefore, there exists $w\in H_0^1(\Omega)\setminus\{0\}$
	such that
	\begin{equation}\label{min levelsN-effective-test}
		\frac{
			\displaystyle
			\int_\Omega
			\left(
			|\nabla w|^2+\lambda_{\bs a}w^2
			\right)\,dx
		}{
			\displaystyle
			\left(
			\int_\Omega|w|^{2^*}\,dx
			\right)^{2/2^*}
		}
		<S.
	\end{equation}
	
	We now consider	$\mf w_{\bs a}
	:=
	(a_1w,\dots,a_Kw)\in H$. Then
	\begin{equation*}
	Q(\mf w_{\bs a})
	:=
	\sum_{i=1}^K
	\int_\Omega
	\left(
	|\nabla(a_iw)|^2+\lambda_i(a_iw)^2
	\right)\,dx =
	|\bs a|^2
	\int_\Omega
	\left(
	|\nabla w|^2+\lambda_{\bs a}w^2
	\right)\,dx,
	\end{equation*}
	whereas
	\begin{equation*}
	P(\mf w_{\bs a})
	:=
	\int_\Omega
	\left(
	\sum_{i=1}^K\mu_i|a_iw|^{2^*}
	+
	K\beta\prod_{i=1}^K|a_iw|^q
	\right)\,dx
	=
	\left(\sum_{i=1}^K\mu_i a_i^{2^*}
	+
	K\beta\prod_{i=1}^K a_i^q\right)
	\int_\Omega|w|^{2^*}\,dx.
	\end{equation*}
	
	On the other hand, since $\mf U=(a_1U_{1,0},\dots,a_KU_{1,0})\in\cN'$,
	we have
	\begin{multline*}
	    \sum_{i=1}^K a_i^2 = S^{-N/2}\left(\sum_{i=1}^K a_i^2\right) \int_{\R^N}|\nabla U_{1,0}|^2\,dx\\
        =S^{-N/2}\left(\sum_{i=1}^K\mu_i a_i^{2^*}	+	K\beta\prod_{i=1}^K a_i^q\right)\int_{\R^N}|U_{1,0}|^{2^*}\,dx = \sum_{i=1}^K\mu_i a_i^{2^*}+	K\beta\prod_{i=1}^K a_i^q
	\end{multline*}
	
	Furthermore, since $\mf U$ is a ground state of
	\eqref{eq:system D12},
	\begin{equation*}
		\cA'
		=
		\frac1N
		S^{N/2} |\bs a|^2.
	\end{equation*}
	
	Let $t>0$ be the unique number such that $t\mf w_{\bs a}\in\mathcal M'$. Then $t^2Q(\mf w_{\bs a})
	=
	t^{2^*}P(\mf w_{\bs a}),
	$ and hence, by the Nehari identity,
	\begin{equation*}
	E(t\mf w_{\bs a})
	=
	\frac1N
	\frac{
		Q(\mf w_{\bs a})^{\frac{2^*}{2^*-2}}
	}{
		P(\mf w_{\bs a})^{\frac{2}{2^*-2}}
	}.
	\end{equation*}
	Therefore, recalling \eqref{min levelsN-effective-test}, we obtain
	\begin{align*}
	\cB' \leq 
	\frac1N
	\frac{
		\left(
		|\bs a|^2
		\displaystyle
		\int_\Omega
		\left(
		|\nabla w|^2+\lambda_{\bs a}w^2
		\right)dx
		\right)^{N/2}
	}{
		\left(
		| \bs a|^2
		\displaystyle
		\int_\Omega|w|^{2^*}\,dx
		\right)^{(N-2)/2}
	}=
	\frac{|\bs a|^2}{N}
	\left[
		\frac{
			\displaystyle
			\int_\Omega
			\left(
			|\nabla w|^2+\lambda_{\bs a}w^2
			\right)dx
		}{
			\displaystyle
			\left(
			\int_\Omega|w|^{2^*}\,dx
			\right)^{2/2^*}
		}
	\right]^{N/2} <
	\frac{|\bs a|^2}{N}S^{N/2}
	=
	\cA',
	\end{align*}
	which concludes the proof.
\end{proof}

We will make use of the following Brezis-Lieb-type lemma.
\begin{lemma}\label{lem:BL-interaction}
Let $q\geq1$ with $Kq=2^*$, and let
$\mf u_n,\mf u\in H$ be such that
$\mf u_n\rightharpoonup\mf u$ weakly in $H$.
Setting $\mf w_n:=\mf u_n-\mf u$, we have, up to a subsequence,
\[
\int_\Omega
\left(
\prod_{i=1}^K|u_{i,n}|^q
-
\prod_{i=1}^K|w_{i,n}|^q
-
\prod_{i=1}^K|u_i|^q
\right)\,dx
\longrightarrow0.
\]
Analogously, up to a subsequence,
\[
\int_\Omega
\left(
\prod_{i=1}^K u_{i,n}
-
\prod_{i=1}^K w_{i,n}
-
\prod_{i=1}^K u_i
\right)\,dx
\longrightarrow0.
\]
\end{lemma}

\begin{proof}
We prove the first result, as the second is analogous. Set
\[
F(\mf z):=\prod_{i=1}^K|z_i|^q,
\qquad \mf z\in\R^K.
\]
We first claim that, for every $\delta>0$,
there exists $C_\delta>0$ such that
\begin{equation}\label{eq:BL-product-estimate}
|F(\mf a+\mf b)-F(\mf a)|
\leq
\delta|\mf a|^{2^*}
+C_\delta|\mf b|^{2^*}
\qquad
\forall\mf a,\mf b\in\R^K.
\end{equation}

Indeed, fix $\delta>0$. The homogeneity of $F$
and its uniform continuity on bounded sets yield
\[
|F(\mf a+\mf b)-F(\mf a)|=|\mf a|^{2^*}\left|F\left(\frac{\mf a+\mf b}{|\mf a|}\right)-F\left(\frac{\mf a}{|\mf a|}\right)\right|
\leq \delta|\mf a|^{2^*},
\]
if $|\mf b|\leq\eta|\mf a|$ for some $\eta=\eta(\delta)$.
On the other hand, if $|\mf b|>\eta|\mf a|$,
using $|F(\mf z)|\leq C|\mf z|^{2^*}$, we obtain
\[
|F(\mf a+\mf b)-F(\mf a)|
\leq C\bigl(|\mf a+\mf b|^{2^*}
+|\mf a|^{2^*}\bigr)
\leq C_\delta|\mf b|^{2^*},
\]
and this proves the claim.

Up to a subsequence, $\mf w_n\to\bs0$
almost everywhere in $\Omega$.
Applying \eqref{eq:BL-product-estimate}
with $\mf a=\mf w_n$ and $\mf b=\mf u$,
we obtain
\[
|F(\mf u_n)-F(\mf w_n)-F(\mf u)|\leq
\delta|\mf w_n|^{2^*}
+C_\delta|\mf u|^{2^*}
+|F(\mf u)|.
\]
Consequently, the functions
\[
\left(
|F(\mf u_n)-F(\mf w_n)-F(\mf u)|
-\delta|\mf w_n|^{2^*}
\right)^+
\]
converge to zero almost everywhere
and are dominated by
$C_\delta|\mf u|^{2^*}+|F(\mf u)|\in L^1(\Omega)$.
By dominated convergence,
\[
\limsup_{n\to\infty}
\int_\Omega
|F(\mf u_n)-F(\mf w_n)-F(\mf u)|\,dx
\leq
\delta\sup_n\|\mf w_n\|_{L^{2^*}}^{2^*}.
\]
Since $\delta>0$ is arbitrary, the
conclusion follows.
\end{proof}

We are now in position to conclude the proof.
\begin{proof}[Proof of Theorem \ref{thm:ground-state-bounded-cooperative}]
		For notational convenience, set
		\[
		\mathcal E
		:=
		\begin{cases}
			E, & q>1,\\
			\widetilde E, & q=1,
		\end{cases}
		\qquad
		\mathfrak M'
		:=
		\begin{cases}
			\mathcal M', & q>1,\\
			\widetilde{\mathcal M}', & q=1,
		\end{cases}
		\]
        where
        \begin{equation}\label{eq:E tilde def}
\tilde E(\mathbf u)
=\frac12\sum_{i=1}^K\|u_i\|_{\lambda_i}^2
-\frac1{2^*}\sum_{i=1}^K\mu_i\int_\Omega|u_i|^{2^*}dx
-\beta \int_\Omega\prod_{i=1}^K u_i dx.
\end{equation}
and
\begin{equation}\label{eq:M tilde prime def}
    \tilde \cM'
=\left\{\mathbf u\in H\setminus\{\mathbf0\}:
\sum_{i=1}^K\|u_i\|_{\lambda_i}^2
=\sum_{i=1}^K\mu_i\int_\Omega|u_i|^{2^*}dx
+K\beta \int_\Omega\prod_{i=1}^K u_i dx\right\}
\end{equation}
The functional $\mathcal E$ is of class $C^1$ on $H$ in both cases. We claim that
\begin{equation}\label{eq4.13}
\cB '>0.
\end{equation}
Indeed, for any $\mf u\in \mathfrak M'$,
the Sobolev embedding and H\"{o}lder's inequality with $Kq=2^*$ shows that 
\[
	\displaystyle\sum_{i=1}^K\mu_i\int_\Omega |u_{i}|^{2^*}dx+K\beta\int_\Omega\prod_{i=1}^K|u_{i}|^qdx
    \leq C\left(\displaystyle\sum_{i=1}^K\norm{u}^2_{\lambda_i}\right)^\frac{2^*}{2},
	\]
which together with the definition of $\mathfrak M'$ yields $$\displaystyle\sum_{i=1}^K\norm{u}^2_{\lambda_i}\geq C>0,
$$
and we conclude that~\eqref{eq4.13} holds.

It is standard to check that $\mathcal E$ has a mountain pass geometry, at the mountain pass level $\cB'>0$;
hence, the Mountain Pass Theorem (see \cite[Theorem 2.9, 2.10]{Wil96}) implies that 
there exists $\{\mf u_n\}\subset\mathfrak M'$ such that 
\begin{equation*}
	\mathcal E(\mf u_n)\to \cB ', \quad \mathcal E'(\mf u_n)\to 0
\end{equation*}
as $n\to\infty$.

Therefore $\mf u_n$ is bounded and thus, up to a subsequence,
	\[
	u_{i,n}\rightharpoonup u_i
	\quad\text{weakly in }H_0^1(\Omega),
	\]
	\[
	u_{i,n}\to u_i
	\quad\text{strongly in }L^s(\Omega)
	\quad\text{for every }1\leq s<2^\ast,
	\]
    and as a consequence $\mathcal E'(\mf u)=0$. Set	$\mf w_n:=\mf u_n-\mf u$. By the Brezis--Lieb lemma, for every $i=1,\ldots,K$,
	\begin{equation}\label{min levelsL-components}
		\int_\Omega |u_{i,n}|^{2^*}dx
		=
		\int_\Omega |u_i|^{2^*}dx
		+
		\int_\Omega |w_{i,n}|^{2^*}dx
		+o(1).
	\end{equation}
	If $q>1$, Lemma~\ref{lem:BL-interaction} yields
	\[
	\int_\Omega\prod_{i=1}^K|u_{i,n}|^qdx
	=
	\int_\Omega\prod_{i=1}^K|u_i|^qdx
	+
	\int_\Omega\prod_{i=1}^K|w_{i,n}|^qdx
	+o(1).
	\]
	If $q=1$, the same Lemma gives instead
	\[
	\int_\Omega\prod_{i=1}^K u_{i,n}dx
	=
	\int_\Omega\prod_{i=1}^K u_idx
	+
	\int_\Omega\prod_{i=1}^K w_{i,n}dx
	+o(1).
	\]
	Since $\mathcal E'(\mf u_n)\to0$ in $H'$ and
	$\mathcal E'(\mf u)=0$, combining these facts with
	\eqref{min levelsL-components}, we obtain,
	for every $i=1,\ldots,K$,
	\begin{equation}\label{eq:remainder-component}
		\int_\Omega |\nabla w_{i,n}|^2dx
		-
		\int_\Omega
		\left(
		\mu_i|w_{i,n}|^{2^*}
		+
		\beta\prod_{j=1}^K|w_{j,n}|^q
		\right)dx
		=o(1).
	\end{equation}
	Furthermore,
	\begin{equation}\label{eq:energy-splitting}
		\mathcal E(\mf u_n)
		=
		\mathcal E(\mf u)+I(\mf w_n)+o(1),
	\end{equation}
	where
	\[
	I(\mf w)
	:=
	\frac12\sum_{i=1}^K\int_\Omega|\nabla w_i|^2dx
	-
	\frac1{2^*}\sum_{i=1}^K\mu_i\int_\Omega|w_i|^{2^*}dx
	-
	\frac{\beta}{q}\int_\Omega\prod_{i=1}^K|w_i|^qdx.
	\]
	(As above, when $q=1$ all the products $\prod_{j=1}^{K}|w_{j,n}|^{q}$ are understood without absolute values, i.e. as
$\prod_{j=1}^{K} w_{j,n}$, consistently with the definition of $\widetilde{E}$; accordingly, $I$ is replaced by its
smooth counterpart $\widetilde{I}$, and $\mathcal{A}'=\inf_{\widetilde{\mathcal{N}}'}\widetilde{I}$).
	Passing to a subsequence, we may assume that, for some $b_i\geq 0$,
	\[
	\int_\Omega|\nabla w_{i,n}|^2dx\rightarrow b_i,
	\qquad\text{as }n\to\infty.
	\]
	Summing \eqref{eq:remainder-component} over $i=1,\ldots,K$, we obtain
	\[
	I(\mf w_n)
	=
	\frac1N\sum_{i=1}^K b_i+o(1).
	\]
	Letting $n\to+\infty$ in \eqref{eq:energy-splitting}, we conclude that
	\begin{equation}\label{eq:energy-limit-splitting}
		0\leq \mathcal E(\mf u)
		\leq
		\mathcal E(\mf u)+\frac1N\sum_{i=1}^K b_i
		=
		\lim_{n\to+\infty}\mathcal E(\mf u_n)
		=
		\cB '.
	\end{equation}
		
	It remains to show that $\mf u\neq\bs0$. If $\mf u=\bs0$, since $\cB '>0$, \eqref{eq:energy-limit-splitting} implies that $\sum_{i=1}^K b_i>0$. Then, for $n$ large we may assume $\mf w_n\neq \bs 0$. Recalling \eqref{eq:remainder-component}, it is straightforward to see that there exists $t_n>0$ such that $t_n \mf w_n\in \cN'$, and $t_n\to 1$ as $n\to \infty$. Then, by \eqref{eq:energy-limit-splitting}, we have
	$$
	\cB '=\frac{1}{N}\sum_{i=1}^K b_i=\lim_{n\to+\infty} I(\mf w_n)=\lim_{n\to+\infty} I(t_n\mf w_n)\geq \cA',
	$$
	contradicting \eqref{min levelsprime-less-Aprime}. Thus $\mf u\neq\bs0$ and
	$\mf u\in\mathfrak M'$.
	
		Consequently, $\cB'$ is attained by a critical point $\mf u$ of
$\mathcal E$. If $q>1$, recalling that $\cM'$ is the classical 
Nehari manifold and hence a natural constraint, we deduce that 
$\mf u$ is a ground state solution of \eqref{eq:system 1}. Moreover, 
since we may assume without loss of generality that the functions 
in the sequence $\{\mf u_n\}$ are nonnegative, $\mf u$ is nonnegative.

If $q=1$, then $\mf u$ is a ground state solution of
\[
\begin{cases}
-\Delta u_i +\lambda_i u_i
=
\mu_i |u_i|^{2^\ast-2}u_i
+
\beta
\displaystyle\prod_{\substack{j=1\\j\neq i}}^K u_j,
& \text{in }\Omega,\quad i=1,\dots,K,\\[0.4em]
u_i\in H_0^1(\Omega),
& i=1,\dots,K.
\end{cases}
\]
We claim that $|\mf u|=(|u_1|,\dots,|u_K|)$ is also a ground state 
solution. To this end, it suffices to show that 
$|\mf u|\in \mathfrak M'$. Let $t>0$ be such that 
$t|\mf u|\in \mathfrak M'$. Since $\beta>0$ and 
$\mf u\in \mathfrak M'$, it follows from 
\eqref{eq:M tilde prime def} that $t\leq 1$. On the other hand,
\[
\cB'\leq \cE(t|\mf u|)=t^2 E(\mf u)\leq E(\mf u)=\cB',
\]
whence $t=1$, proving the claim. Thus, also in this case, we obtain 
a nonnegative ground state solution.

By elliptic regularity, $\mf u$ is a classical solution. Hence, 
by the strong maximum principle, every nontrivial component of 
$\mf u$ is strictly positive in $\Omega$. In particular, $\mf u$ 
is a solution of \eqref{eq:system 1}. The characterization then 
follows from Proposition~\ref{prop:beta bar charact H01}.
\end{proof}

\section{Least energy positive solutions in the weak cooperative case}\label{sec:weak cooperative}
The aim of this section is to prove Theorems \ref{thm:weak cooperative D12} and \ref{thm:weak cooperative bounded domain} in the first and second subsections, respectively.

\subsection{The weak cooperative case on $\R^N$}
We first introduce an algebraic problem that will be useful in proving Theorem \ref{thm:weak cooperative bounded domain}. For $\beta>0$, we define the set of algebraic constraints
\[
	\Gamma_\beta
	:=
	\left\{
	\mf y=(y_1,y_2,y_3)\in\R_{>0}^3:
	\mu_i y_i^4 + \beta y_j^2 y_k^2 \geq S
	\quad\text{for every }
	\{i,j,k\}=\{1,2,3\}
	\right\}
\]
and the corresponding algebraic minimization problem
\[
	m_\beta
	:=
	\inf_{\mf y\in\Gamma_\beta}
	\sum_{i=1}^3 y_i^2.
\]

\begin{proposition}\label{P5.1}
For $\beta>0$ sufficiently small, $m_\beta$ is attained at some $\mf y \in \R_{>0}^3$. Moreover, all algebraic constraints are active at $\mf y$, that is,
\begin{equation}\label{eq:all-algebraic-constraints-active}
	\mu_i y_i^4 + \beta y_j^2 y_k^2 = S,
	\qquad \{i,j,k\}=\{1,2,3\}.
\end{equation}
\end{proposition}

\begin{proof}
Let $(\mf y_n)\subset\Gamma_\beta$ be a minimizing sequence for $m_\beta$. Testing the constraint with the admissible point
\[
	\left(
	\left(\frac{S}{\mu_1}\right)^{1/4},
	\left(\frac{S}{\mu_2}\right)^{1/4},
	\left(\frac{S}{\mu_3}\right)^{1/4}
	\right)
	\in\Gamma_\beta,
\]
we immediately obtain the upper bound
\[
	m_\beta \leq C_0 := \sqrt{S}\sum_{i=1}^3 \mu_i^{-1/2}.
\]
Consequently, discarding finitely many terms, if necessary, there exists $C>0$ (independent of $\beta$ in a bounded interval) such that
\[
	0 < y_{i,n} \leq C, \qquad i=1,2,3.
\]
To show that the sequence is bounded away from zero, we use the constraint definition:
\[
	S \leq \mu_i y_{i,n}^4 + \beta y_{j,n}^2 y_{k,n}^2 \leq \mu_i y_{i,n}^4 + \beta C^4.
\]
Choosing $\beta_1>0$ small enough such that $\beta_1 C^4 \leq \frac{S}{2}$, it follows that for every $\beta\in(0,\beta_1)$,
\[
	y_{i,n}^4 \geq \frac{S}{2\mu_i},
\]
which yields a uniform lower bound
\begin{equation}\label{eq:yi-lower-bound}
	y_{i,n} \geq c > 0, \qquad i=1,2,3,
\end{equation}
for some $c>0$ independent of $n$ and $\beta\in(0,\beta_1)$. 
Passing to a subsequence, $\mf y_n \to \mf y \in [c,C]^3$. Since $\Gamma_\beta$ is closed, $\mf y \in \Gamma_\beta$ and $\sum_{i=1}^3 y_i^2 = m_\beta$.

\medskip
To complete the proof, it remains to show that all algebraic constraints are active at $\mf y$. We define the constraint functions
\[
	g_i(\mf y) := S - \mu_i y_i^4 - \beta y_j^2 y_k^2, \qquad \{i,j,k\}=\{1,2,3\},
\]
so that $\Gamma_\beta$ is locally characterized by $g_i(\mf y)\leq 0$ for $i=1,2,3$. Since $\mf y \in \R_{>0}^3$ lies strictly in the interior of the positive orthant, the Kuhn--Tucker optimality conditions yield multipliers $\lambda_1,\lambda_2,\lambda_3 \geq 0$ such that
\begin{equation}\label{eq:KKT-stationarity}
	2 y_i = 4\mu_i y_i^3 \lambda_i + 2\beta y_i (y_k^2 \lambda_j + y_j^2 \lambda_k), \qquad \{i,j,k\}=\{1,2,3\},
\end{equation}
together with the complementary slackness conditions
\[
	\lambda_i g_i(\mf y) = 0, \qquad i=1,2,3.
\]

We claim that $\lambda_i > 0$ for every $i=1,2,3$. Suppose by contradiction that $\lambda_1 = 0$. Then, \eqref{eq:KKT-stationarity} implies that
\[
	1 = \beta (y_3^2 \lambda_2 + y_2^2 \lambda_3).
\]
On the other hand, non-negativity of all terms in \eqref{eq:KKT-stationarity} for $i=2,3$ implies
\[
	\lambda_2 \leq \frac{1}{2\mu_2 y_2^2}, \qquad \lambda_3 \leq \frac{1}{2\mu_3 y_3^2}.
\]
Substituting these bounds yields
\[
	1 \leq \frac{\beta}{2} \left( \frac{y_3^2}{\mu_2 y_2^2} + \frac{y_2^2}{\mu_3 y_3^2} \right) \leq C_1 \beta,
\]
where $C_1>0$ is independent of $\beta\in(0,\beta_1)$ thanks to \eqref{eq:yi-lower-bound} and $y_i\leq C$. For $\beta>0$ sufficiently small, this is a contradiction. The same reasoning applies to $\lambda_2$ and $\lambda_3$, proving that $\lambda_i > 0$ for all $i=1,2,3$. Therefore, $g_i(\mf y) = 0$ for $i=1,2,3$, and the proof is complete.
\end{proof}

We can now prove Theorem \ref{thm:weak cooperative D12}.
\begin{proof}[Proof of Theorem \ref{thm:weak cooperative D12}]
Let $\mf u=(u_1,u_2,u_3)\in\cN$ and define $y_i := |u_i|_6 > 0$. By the sharp Sobolev inequality, the component-wise Nehari identities, and H\"older's inequality, we have
\[
	S y_i^2 \leq \int_{\R^3}|\nabla u_i|^2 dx = \mu_i |u_i|_6^6 + \beta \int_{\R^3} |u_1 u_2 u_3|^2 dx \leq \mu_i y_i^6 + \beta y_1^2 y_2^2 y_3^2.
\]
Dividing by $y_i^2 > 0$ yields
\[
	S \leq \mu_i y_i^4 + \beta y_j^2 y_k^2, \qquad \{i,j,k\}=\{1,2,3\}.
\]
Thus $\mf y = (y_1,y_2,y_3) \in \Gamma_\beta$. As a consequence,
\[
	I(\mf u) = \frac{1}{3} \sum_{i=1}^3 \int_{\R^3} |\nabla u_i|^2 dx \geq \frac{S}{3} \sum_{i=1}^3 y_i^2 \geq \frac{S}{3} m_\beta.
\]
Taking the infimum over $\cN$, we conclude that
\begin{equation}\label{eq:A-lower-algebraic}
	\cA \geq \frac{S}{3} m_\beta.
\end{equation}

We now show that equality holds in the preceding inequality and that $\cA$ is attained. Let $\mf y=(y_1,y_2,y_3)$ be the minimizer of $m_\beta$ obtained in Proposition \ref{P5.1}. Moreover,
\[
	S = \mu_i y_i^4 + \beta y_j^2 y_k^2, \qquad i=1,2,3.
\]
Let $U=U_{\varepsilon,y}$ be an Aubin--Talenti bubble normalized such that
\[
	\int_{\R^3} |\nabla U|^2 dx = \int_{\R^3} U^6 dx = S^{3/2}.
\]
We set
\[
	u_i := \frac{y_i}{S^{1/4}} U.
\]
Direct calculations show that
\[
	|u_i|_6^2 = \frac{y_i^2}{S^{1/2}} S^{1/2} = y_i^2, \qquad \int_{\R^3} |\nabla u_i|^2 dx = \frac{y_i^2}{S^{1/2}} S^{3/2} = S y_i^2,
\]
and
\[
	\int_{\R^3} |u_i|^6 dx = \frac{y_i^6}{S^{3/2}} S^{3/2} = y_i^6, \qquad \int_{\R^3} |u_1 u_2 u_3|^2 dx = \frac{y_1^2 y_2^2 y_3^2}{S^{3/2}} S^{3/2} = y_1^2 y_2^2 y_3^2.
\]
By \eqref{eq:all-algebraic-constraints-active}, $\mf u \in \cN$. Furthermore,
\[
	I(\mf u) = \frac{1}{3} \sum_{i=1}^3 \int_{\R^3} |\nabla u_i|^2 dx = \frac{S}{3} \sum_{i=1}^3 y_i^2 = \frac{S}{3} m_\beta.
\]
Combining this with \eqref{eq:A-lower-algebraic}, we conclude that $\cA = \frac{S}{3} m_\beta$ and that $\cA$ is achieved at $\mf u$. The proof that $\mf u$ is a critical point (and hence a solution) follows by arguing as in \cite[Proposition 4.7]{Giaretto_Soave_2026}.    
\end{proof}

\subsection{The weak cooperative case on bounded domains}
Recalling the definition of $\cB$ in \eqref{min levels}, since
\[
\int_\Omega\left(|\nabla u|^2+\lambda_i u^2\right)dx
\geq
\left(1+\frac{\lambda_i}{\lambda_1(\Omega)}\right)
\int_\Omega|\nabla u|^2dx,
\qquad i=1,\dots,K,
\]
it is standard to see that $\cB>0$.

The strategy is to minimize on $\cM$ while ruling out concentration phenomena. The first step is to establish the following energy estimate:
\begin{equation}\label{eq:upper bound B}\tag{E}
	\cB<\cA \quad \text{and} \quad \cB<\min_{i=1,\dots,K} \left\{
		\sum_{j\neq i} B_{\mu_j}
		+\frac1N\mu_i^{-\frac{N-2}{2}}S^{\frac N2}
		\right\},
\end{equation}
where $B_\mu$ is defined in \eqref{eq:B_mu def least energy bounded domain}.

\begin{lemma}\label{lem:upper-bound-B-N3-beta-positive}
	Let $N=3$, $K=3$, $q=2$, and assume that \eqref{coop-intro:linear-assumptions} holds. If $\beta>0$ is sufficiently small, then \eqref{eq:upper bound B} holds.
\end{lemma}

\begin{proof}
We first prove 
\begin{equation}\label{eq:upper bound 1}
	\cB<\sum_{i=1}^3 B_{\mu_i}.
\end{equation}
For every $i=1,2,3$, let $z_i$ be a least energy positive solution of
\[
\begin{cases}
	-\Delta z_i+\lambda_i z_i = \mu_i z_i^5 &\text{in }\Omega,\\
	z_i=0 &\text{on }\partial\Omega.
\end{cases}
\]
Since $z_i$ belongs to the scalar Nehari manifold,
\[
	\|z_i\|_{\lambda_i}^2
	=
	\mu_i\int_\Omega z_i^6\,dx, \quad \text{and} \quad B_{\mu_i}
	=
	\frac13\|z_i\|_{\lambda_i}^2,
	\qquad i=1,2,3.
\]

Consider the fiber map
\[
	\Psi(\bs t):=E(t_1z_1, t_2 z_2, t_3 z_3),
	\qquad
	\bs t=(t_1,t_2,t_3)\in[0,+\infty)^3.
\]
Then
\begin{equation*}
	\Psi(\bs t)
	=
	\sum_{i=1}^3
	\|z_i\|_{\lambda_i}^2\left(\frac{t_i^2}{2}-\frac{t_i^6}{6}\right)
	-\frac{\beta}{2}
	t_1^2 t_2^2 t_3^2 \int_\Omega z_1^2 z_2^2 z_3^2 \,dx.
\end{equation*}
For every $t_i\geq0$,
\[
	\frac{t_i^2}{2}-\frac{t_i^6}{6}\leq \frac{1}{3},
\]
and equality holds if and only if $t_i=1$. Therefore, for every $\bs t\in(0,+\infty)^3$,
\begin{equation*}
	\Psi(\bs t)<\frac{1}{3}\sum_{i=1}^3 \|z_i\|_{\lambda_i}^2 = \sum_{i=1}^3 B_{\mu_i}.
\end{equation*}

We now choose $\bs t=(t_1,t_2,t_3)\in(0,+\infty)^3$ such that
\[
	(t_1z_1,t_2z_2,t_3z_3)\in\mathcal M.
\]
To see that such a $\bs t$ exists, consider the functions $F_i : (0,+\infty)^3 \to \R$ defined for $i=1,2,3$ by
\[
	F_i(\bs t) := t_i^2 \|z_i\|_{\lambda_i}^2 - \mu_i t_i^6 \int_\Omega z_i^6 \,dx - \beta t_1^2 t_2^2 t_3^2 \int_\Omega z_1^2 z_2^2 z_3^2 \,dx.
\]
Recalling that $\|z_i\|_{\lambda_i}^2 = \mu_i \int_\Omega z_i^6 \,dx$, we can divide by $t_i^2$ to obtain
\[
	\frac{F_i(\bs t)}{t_i^2} = \|z_i\|_{\lambda_i}^2 \left(1 - t_i^4\right) - \beta \left(\prod_{j\neq i} t_j^2\right) \int_\Omega z_1^2 z_2^2 z_3^2 \,dx.
\]
For $r>0$ sufficiently small and $R>1$ sufficiently large, independently of the other components $t_j \in [r,R]$ ($j \neq i$), we have $F_i(\bs t)>0$ when $t_i=r$ and $F_i(\bs t)<0$ when $t_i=R$. Thus, Miranda's Theorem yields a vector $\bs t \in (r,R)^3$ such that $F_i(\bs t)=0$ for $i=1,2,3$, which precisely means that $(t_1z_1,t_2z_2,t_3z_3)\in\mathcal M$.

Then, by the definition of $\cB$,
\[
	\cB
	\leq
	E(t_1z_1,t_2 z_2,t_3 z_3)
	= \Psi(\bs t)
	<
	\sum_{i=1}^3 B_{\mu_i}.
\]
This proves \eqref{eq:upper bound 1}.

\medskip

We now show that $\cB<\cA$.
Recall that, by Theorem \ref{thm:weak cooperative D12}, for every $\beta>0$ sufficiently small $\cA$ is attained by a synchronized configuration $(\ell_1(\beta)U,\ell_2(\beta)U,\ell_3(\beta)U)$, where $U$ is an Aubin--Talenti function normalized so that
\[
	\int_{\R^3}|\nabla U|^2dx
	=
	\int_{\R^3}U^6dx
	=
	S^{3/2},
\]
$\ell_i(\beta)>0$ is bounded uniformly in $\beta$ and satisfies
\begin{equation}\label{eq:ell-system-A}
	\ell_i(\beta)^2
	=
	\mu_i\ell_i(\beta)^6
	+
	\beta\ell_1(\beta)^2 \ell_2(\beta)^2 \ell_3(\beta)^2,
	\qquad i=1,2,3.
\end{equation}
In particular,
\begin{equation}\label{eq:A-synchronized}
	\cA=\cA(\beta)
	=
	\frac{1}{3} S^{3/2}\sum_{i=1}^3\ell_i(\beta)^2.
\end{equation}

We aim at constructing a suitable family of competitors. First, set $\lambda_{\max}:=\max_{1\le i \le 3}\lambda_i$, then, since $\lambda_{\max}<-\lambda^*(\Omega)$, by the usual construction in \cite{BN} there exists a family of nonnegative functions $w_\varepsilon\in H_0^1(\Omega)$ such that, for $\varepsilon>0$ sufficiently small and for every $i=1,2,3$,
\begin{equation}\label{eq:quotient less than S}
	\frac{
	\displaystyle
	\int_\Omega
	\left(
	|\nabla w_\varepsilon|^2
	+\lambda_iw_\varepsilon^2
	\right)dx}
	{\displaystyle
	\left(
	\int_\Omega|w_\varepsilon|^6
	dx\right)^{1/3}}
	\leq\frac{
	\displaystyle
	\int_\Omega
	\left(
	|\nabla w_\varepsilon|^2
	+\lambda_{\max}w_\varepsilon^2
	\right)dx}
	{\displaystyle
	\left(
	\int_\Omega|w_\varepsilon|^6 dx
	\right)^{1/3}}
	<S,
\end{equation}
and
\begin{equation}\label{eq:quotient tends to S}
	\frac{
	\displaystyle
	\int_\Omega
	\left(
	|\nabla w_\varepsilon|^2
	+\lambda w_\varepsilon^2
	\right)dx}
	{\displaystyle
	\left(
	\int_\Omega|w_\varepsilon|^6 dx
	\right)^{1/3}}
	\longrightarrow S
	\qquad\text{as }\varepsilon\to0
\end{equation}
for every fixed $\lambda\in\R$. We aim at finding $\bs s_\eps=(s_{1,\eps},s_{2,\eps},s_{3,\eps})\in \R^3_{>0}$ such that $\bs s_\eps w_\eps\in \cM$ and $E(\bs s_\eps w_\eps)<\cA$ for $\eps>0$ small enough. 

For $\bs t, \bs a\in\R^3$, consider
\[
	J(\bs t,\bs a)
	:=
	\frac12\sum_{i=1}^3 a_i t_i^2
	-
	\frac{1}{6}\sum_{i=1}^3\mu_i t_i^6
	-
	\frac{\beta}{2} t_1^2 t_2^2 t_3^2.
\]
For $\bs t\in\R_{>0}^3$, the condition $\nabla_{\bs t}J(\bs t,\bs a)=0$ is equivalent to
\begin{equation*}
	a_i
	=
	\mu_i t_i^4
	+
	\beta \prod_{j\neq i} t_j^2,
	\qquad
	i=1,2,3.
\end{equation*}
We apply the Implicit Function Theorem to
\[
	F_{i}(\bs t,\bs a)
	:=
	\mu_i t_i^4+\beta \prod_{j\neq i} t_j^2-a_i,
	\qquad i=1,2,3.
\]
Notice that \eqref{eq:ell-system-A} is equivalent to
\begin{equation*}
	1
	=
	\mu_i \ell_i(\beta)^4
	+
	\beta \prod_{j\neq i} \ell_j(\beta)^2,
	\qquad
	i=1,2,3,
\end{equation*}
therefore
\[
	\bs F(\bs \ell(\beta),\bs a_0)=0,
	\qquad
	\bs a_0:=(1,1,1).
\]
Moreover \eqref{eq:ell-system-A} implies that, as $\beta\to0^+$,
\[
	\ell_i(\beta)\to \mu_i^{-1/4},
	\qquad i=1,2,3.
\]
Hence the Jacobian of $\bs F$ with respect to $\bs t$ at $(\bs \ell(\beta),\bs a_0)$ satisfies
\[
	D_{\bs t}\bs F(\bs \ell(\beta),\bs a_0)
	\longrightarrow
	\operatorname{diag}
	\left(
	4\mu_1^{1/4}, 4\mu_2^{1/4}, 4\mu_3^{1/4}
	\right) \qquad \text{as } \b\to 0.
\]
Thus, reducing $\beta>0$ if necessary, $D_{\bs t}\bs F(\bs \ell(\beta) ,\bs a_0)$ is invertible.

The Implicit Function Theorem therefore yields a neighborhood $\mathcal U$ of $\bs a_0$ and a $C^1$ map $\bs t:\mathcal U\to\R_{>0}^3$ such that $\bs t(\bs a_0)=\bs \ell(\beta)$ and
\begin{equation}\label{eq:a(s) critical point J}
	\nabla_{\bs t}J(\bs t(\bs a),\bs a)=0 \quad \forall \bs a\in \cU.
\end{equation}
Define
\[
	a_{i}(\eps)=\frac{\norm{w_\eps}_{\lambda_i}^2}{S |w_\eps|_6^2}, \quad s_{i,\eps}=\frac{S^{1/4}t_i(\bs a(\eps))}{|w_\eps|_6}.
\]
Notice that, by \eqref{eq:quotient less than S} and \eqref{eq:quotient tends to S}, $a_i(\eps)\to 1$ as $\eps\to 0$ and $a_i(\eps)<1$ for $\eps$ small enough. Moreover, by construction,
\begin{align*}
	s_{i,\eps}^2\norm{w_\eps}_{\lambda_i}^2
	&=
	S^{1/2}t_i(\bs a(\varepsilon))^2 \frac{\|w_{\varepsilon}\|_{\lambda_i}^2}{|w_\eps|_6^2} = S^{3/2} t_i(\bs a(\varepsilon))^2 a_{i}(\eps)\\
	&=
	S^{3/2}\left(\mu_i t_i(\bs a(\varepsilon))^6
	+
	\beta t_1(\bs a(\varepsilon))^2 t_2(\bs a(\varepsilon))^2 t_3(\bs a(\varepsilon))^2 \right)\\
	&=\mu_i s_{i,\eps}^6
	\int_\Omega|w_\eps|^6dx
	+
	\beta s_{1,\eps}^2s_{2,\eps}^2s_{3,\eps}^2
	\int_\Omega
	|w_\eps|^6 dx,
\end{align*}
thus $\bs s_\eps w_\eps\in \cM$. We now compare the corresponding energy levels.
Recalling \eqref{eq:a(s) critical point J}, we have, for every $\bs a\in \cU$,
\begin{equation}
	\frac{\partial J}{\partial a_i}(\bs t(\bs a),\bs a)
	=
	\frac{t_i(\bs a)^2}{2}
	+
	\left\langle
	\nabla_{\bs t}J(\bs t(\bs a), \bs a),
	\frac{\partial\bs t}{\partial a_i}(\bs a)
	\right\rangle= \frac{t_i(\bs a)^2}{2}>0.
	\label{eq:c-derivative}
\end{equation}
Moreover, notice that, by the normalization of $U$ and by \eqref{eq:ell-system-A} and \eqref{eq:A-synchronized},
\[
  I(\boldsymbol{\ell}(\beta)U)=S^{3/2}J(\boldsymbol{\ell}(\beta),\boldsymbol{a}_0)
  =\frac13 S^{3/2}\sum_{i=1}^{3}\ell_i(\beta)^2=\mathcal{A},
  \qquad \boldsymbol{a}_0=(1,1,1).
\]

For $\varepsilon>0$ sufficiently small such that the convex combination of $\bs a_0$ and $\bs a(\eps)$ is contained in $\cU$, using \eqref{eq:c-derivative} and denoting by $\bs \gamma(r):=\bs a_0+r(\bs a(\eps)-\bs a_0)$, $r\in [0,1]$,
	\begin{align*}
	E(\bs s_\eps w_\eps)-I(\bs \ell(\beta) U)&=\frac{1}{2}\sum_{i=1}^3 s_{i,\eps}^2 \norm{w_\eps}_{\lambda_i}^2 - \frac{1}{6} \sum_{i=1}^3\mu_i s_{i,\eps}^6 |w_\eps|_6^6 - \frac{\b}{2}s_{1,\eps}^2 s_{2,\eps}^2 s_{3,\eps}^2 |w_\eps|_6^6- S^{3/2} J(\bs \ell(\b), \bs a_0) \\
    &=S^{3/2}\left[J(\bs t(\bs a(\eps)), \bs a(\eps))-J(\bs \ell(\b), \bs a_0)\right]\\
    &=S^{3/2} \int_0^1 \frac{d}{dr}J(\bs t(\bs \gamma(r)),\bs \gamma(r)) dr\\
    &=S^{3/2} \int_0^1 \left[\nabla_{\bs t} J(\bs t(\bs \gamma(r)),\bs \gamma(r)) \cdot \frac{d\bs t}{dr}(\bs \gamma(r))+\nabla_{\bs a} J(\bs t(\bs \gamma(r)),\bs \gamma(r))\cdot \frac{d\bs \gamma(r)}{dr}\right]\, dr \\
    &=S^{3/2} \int_0^1 \sum_{i=1}^3 \frac{\partial J}{\partial a_i}(\bs t(\bs \gamma(r)),\bs \gamma(r)) (a_{i}(\eps)-1) \, dr\\
    & =S^{3/2} \int_0^1 \sum_{i=1}^3 \frac{t_i(\bs \gamma(r))^2}{2}(a_i(\eps)-1)\,dr<0,
\end{align*}
where in the last inequality we used \eqref{eq:c-derivative} and $a_{i}(\eps)<1$ for $\eps$ small enough. Hence,
\[
	\cB
	\leq
	E(\mf u_\varepsilon)
	<
	I(\bs \ell(\b)U)
	=
	\cA. \qedhere
\]
\end{proof}

We now collect some preliminary results concerning the component-wise
Nehari manifold. They are direct adaptations to the critical case
$Kq=2^*$ of the corresponding results proved in
\cite{Giaretto_Soave_2026}. Since the proofs are identical, we only
recall the main arguments when useful.

For $i=1,\dots,K$, define
\[
G_{\beta,i}(\mf u)
:=\|u_i\|_{\lambda_i}^2
-\mu_i |u_i|_{2^*}^{2^*}
-\beta P(\mf u), \quad \text{where} \quad P(\mf u)=
\begin{cases}
 \int_\Omega\prod_{j=1}^K |u_j|^q\,dx    & \text{if $q>1$}, \\
\int_\Omega\prod_{j=1}^K u_j\,dx & \text{if $q=1$}
\end{cases}
\]
so that
\[
\mathcal M
=
\left\{
\mf u\in H:
u_i\not\equiv0,\ 
G_{\beta,i}(\mf u)=0
\quad\text{for every }i=1,\dots,K
\right\}.
\]
For $\mf u\in\mathcal M$ and $j=1,\dots,K$, let
\[
(\bs\varphi_j(\mf u))_k
=
\begin{cases}
u_j, & k=j,\\
0, & k\neq j,
\end{cases}
\qquad k=1,\dots,K.
\]

We define $\bs M_\beta(\mf u)\in\R^{K\times K}$ by
\[
(\bs M_\beta(\mf u))_{ij}
:=
G_{\beta,i}'(\mf u)[\bs\varphi_j(\mf u)].
\]

Therefore, in both cases $q>1$ and $q=1$, using $Kq=2^*$ and
$G_{\beta,i}(\mf u)=0$, we have
\[
(\bs M_\beta(\mf u))_{ij}
=
\begin{cases}
(2-2^*)\mu_i|u_i|_{2^*}^{2^*}
+(2-q)\beta P(\mf u),
& i=j,\\[0.5em]
-q\beta\displaystyle\int_\Omega
\prod_{k=1}^K|u_k|^q\,dx,
& i\neq j.
\end{cases}
\]

\begin{lemma}\label{lem:uniform-bounds-components}
		Let $q=2$ and $\b>0$ sufficiently small. Then there exist constants $C_1,C_2>0$ such that
		\[
		C_1\leq |u_i|_{2^*}\leq C_2
		\qquad
		\text{for every } \mf u\in\mathcal M
		\text{ such that } E(\mf u)\leq \cA,
		\]
		for every $i=1,\ldots,K$.
	\end{lemma}
	
	\begin{proof}
		The proof is identical to that of \cite[Lemma 4.5]{Giaretto_Soave_2026}.
	\end{proof}

\begin{lemma}\label{lemma:M_beta_negative}
Let $q=2$ and $\b>0$ sufficiently small. Then, there exists $\eps>0$ such that, for every $\mf u\in\mathcal M$ such that $E(\mf u)\leq \cB+\eps$, $\bs M_\beta(\mf u)$ is strictly negative definite.
\end{lemma}

\begin{proof}
The proof is identical to that of
\cite[Lemma~4.3]{Giaretto_Soave_2026}.
\end{proof}

\begin{lemma}\label{lemma:M_manifold}
Let $q=2$ and $\b>0$ sufficiently small. Then $\mathcal M \cap \{E < \mathcal{B}+\eps\}$ is a smooth manifold of
codimension $K$ in $H$.
\end{lemma}

\begin{proof}
The proof is identical to that of
\cite[Lemma~4.4]{Giaretto_Soave_2026}. Indeed,
\Cref{lemma:M_beta_negative} implies that the differential of
\[
\mf u\longmapsto
\bigl(G_{\beta,1}(\mf u),\dots,G_{\beta,K}(\mf u)\bigr)
\]
is surjective at every $\mf u\in\mathcal M$.
\end{proof}

\begin{proposition}\label{prop:minimizers_are_solutions}
Let $q=2$ and $\b>0$ sufficiently small and let $\mf u\in\mathcal M$ be a minimizer
of $E$ on $\mathcal M$. Then $\mf u$ is a critical point of $E$
and hence a weak solution of \eqref{eq:system 1}.
\end{proposition}

\begin{proof}
The proof is identical to that of
\cite[Proposition~4.7]{Giaretto_Soave_2026}. By
\Cref{lemma:M_manifold}, the Lagrange multiplier rule applies,
while \Cref{lemma:M_beta_negative} implies that all the Lagrange
multipliers vanish.
\end{proof}

Combining the energy estimate \eqref{eq:upper bound B} with the previous lemmas we can conclude the proof.

\begin{proof}[Proof of Theorem \ref{thm:weak cooperative bounded domain}]
    From now on, let $\b>0$ be small enough such that all previous lemmas hold true. Moreover, we take $\eps>0$ so small that $\mathcal{B}+\eps<\mathcal{A}$ (this choice is possible, thanks to \eqref{eq:upper bound B}). Recall that $E$ is coercive and bounded from below on $\mathcal M$, and observe that $\mathcal M \cap \{E \le \mathcal{B}+\eps/2\}$ is a complete metric space. By Ekeland's variational principle, there exists a sequence
$(\mf u_n)\subset\mathcal M \cap \{E \le \mathcal{B}+\eps/2\}$ such that
\begin{equation*}
    E(\mf u_n)\to \cB,
    \qquad
    (E|_{\mathcal M})'(\mf u_n)\to0.
\end{equation*}
In particular, $(\mf u_n)$ is bounded in $H$ and, up to a subsequence,
\[
\mf u_n\rightharpoonup\mf u
\qquad\text{in }H.
\]
Using Lemma \ref{lemma:M_beta_negative} and arguing as in \cite[Lemma 4.6]{Giaretto_Soave_2026}, one can
show that $(\mf u_n)$ is actually an unconstrained Palais--Smale sequence, namely
\[
E'(\mf u_n)\to0
\qquad\text{in }H'.
\]
By the weak convergence, the Sobolev embeddings, and standard weak-continuity arguments for the nonlinear terms, we obtain
\[
E'(\mf u)=0.
\]

Set
\[
\mf w_n:=\mf u_n-\mf u.
\]
By the Brezis--Lieb lemma, for every $i=1,\ldots,K$,
\begin{equation*}
    \|u_{i,n}\|_{2^*}^{2^*}=\|u_i\|_{2^*}^{2^*}+\|w_{i,n}\|_{2^*}^{2^*}+o(1).
\end{equation*}
Moreover, Lemma~\ref{lem:BL-interaction} with $q=2$ gives
\[
\int_\Omega
\left(
\prod_{i=1}^K |u_{i,n}|^2
-
\prod_{i=1}^K |w_{i,n}|^2
-
\prod_{i=1}^K |u_i|^2
\right)\,dx
=o(1).
\]
Since $w_{i,n}\to0$ strongly in $L^2(\Omega)$, combining these relations with $E'(\mf u_n)[\mf u_n]\to0$ and $E'(\mf u)=0$ (resp. $\mf g_{\bs u}=0$) yields, for every $i=1,\ldots,K$,
\begin{equation}\label{eq:remainder-nehari}
    \int_\Omega |\nabla w_{i,n}|^2dx
    -
    \mu_i\int_\Omega |w_{i,n}|^{2^*}dx
    -
    \beta\int_\Omega\prod_{\ell=1}^K|w_{\ell,n}|^2dx
    =o(1).
\end{equation}
Furthermore,
\begin{equation}\label{eq:energy-splitting2}
    E(\mf u_n)
    =
    E(\mf u)+I(\mf w_n)+o(1),
\end{equation}
where, recalling $2^* = 2K$,
\[
I(\mf w)
:=
\frac12\sum_{i=1}^K\int_\Omega|\nabla w_i|^2dx
-
\frac1{2^*}\int_\Omega
\left(
\sum_{i=1}^K\mu_i|w_i|^{2^*}
+
K\beta\prod_{i=1}^K|w_i|^2
\right)dx.
\]

Up to a further subsequence, we may assume that
\[
\int_\Omega|\nabla w_{i,n}|^2dx\to b_i\geq0,
\qquad i=1,\ldots,K.
\]
Summing \eqref{eq:remainder-nehari} over $i$ and recalling that $2K=2^*$, we obtain
\[
I(\mf w_n)
=
\frac1N\sum_{i=1}^K b_i+o(1).
\]
Since $E'(\mf u)=0$, we also have $E(\mf u)\geq0$. Hence, letting $n\to\infty$ in \eqref{eq:energy-splitting2},
\begin{equation}\label{min levels-splitting}
    0\leq E(\mf u)
    \leq
    E(\mf u)+\frac1N\sum_{i=1}^K b_i
    =
    \cB.
\end{equation}

It remains to prove that
\[
u_i\not\equiv0
\qquad\text{for every }i=1,\ldots,K.
\]

\medskip

\noindent\textbf{Case 1: $\mf u=\bs0$.}

\medskip

In this case $w_{i,n}=u_{i,n}$ and $b_i>0$ for every $i$. Indeed, if $b_i=0$ for some $i$, then $u_{i,n}\to0$ in $H_0^1(\Omega)$, contradicting Lemma \ref{lem:uniform-bounds-components}. We claim that there exist $t_{i,n}>0$, $i=1,\ldots,K$, such that
\[
(t_{1,n}w_{1,n},\ldots,t_{K,n}w_{K,n})\in\mathcal N.
\]
For $\bs t=(t_1,\ldots,t_K)\in(0,+\infty)^K$, define
\[
F_{j,n}(\bs t)
:=
t_j^2\int_\Omega|\nabla w_{j,n}|^2dx
-
t_j^{2^*}\mu_j\int_\Omega|w_{j,n}|^{2^*}dx
-
\beta t_j^2\prod_{\ell\neq j}t_\ell^2
\int_\Omega\prod_{\ell=1}^K|w_{\ell,n}|^2dx.
\]
Since $(\mf w_n)$ is bounded, after passing to a subsequence we may assume that, for every $j=1,\ldots,K$,
\[
\mu_j\int_\Omega|w_{j,n}|^{2^*}dx\to c_j>0,
\qquad
\beta\int_\Omega\prod_{\ell=1}^K|w_{\ell,n}|^2dx\to d\geq0.
\]
By \eqref{eq:remainder-nehari},
\begin{equation}\label{min levelsi-ci-d}
    c_j+d=b_j.
\end{equation}

We can choose $0<r<1<R$ so that, for $n$ sufficiently large and uniformly with respect to the remaining variables in $[r,R]$ and with respect to $\b$,
\[
F_{j,n}(\bs t)>0
\quad\text{if }t_j=r,
\qquad
F_{j,n}(\bs t)<0
\quad\text{if }t_j=R.
\]
Miranda's theorem therefore yields $\bs t_n=(t_{1,n},\ldots,t_{K,n})\in[r,R]^K$ such that
\[
F_{j,n}(\bs t_n)=0,
\qquad j=1,\ldots,K.
\]
Up to a subsequence, $t_{i,n}\to\bar t_i\in[r,R]$.
Passing to the limit gives
\begin{equation}\label{eq:t-limit-system}
    \bar t_i^2 b_i
    =
    \bar t_i^{2^*}c_i+d\prod_{\ell=1}^K\bar t_\ell^2,
    \qquad i=1,\ldots,K.
\end{equation}
In particular, dividing by $\bar t_i^2$,
\begin{equation*}
    \bar t_i
    \leq
    \left(\frac{b_i}{c_i}\right)^{\frac1{2^*-2}}.
\end{equation*}

We claim that $\bar t_i=1$ for every $i=1,\ldots,K$. If $d=0$, this follows immediately from \eqref{min levelsi-ci-d} and \eqref{eq:t-limit-system}. Suppose therefore $d>0$. Dividing \eqref{eq:t-limit-system} by $\bar t_i^2$ and using \eqref{min levelsi-ci-d}, we rewrite the limit system for $\bar{\bs t} = (\bar t_1,\ldots,\bar t_K) \in [r,R]^K$ as
\begin{equation}\label{eq:fixed-point-system}
    c_i \left(1 - \bar t_i^{2^*-2}\right) = d \left(\prod_{\ell \neq i}\bar t_\ell^2 - 1\right), \qquad i=1,\ldots,K.
\end{equation}
Clearly, $\bar t_1 = \cdots = \bar t_K = 1$ is a solution of \eqref{eq:fixed-point-system}. We now show that it is the unique solution for $\beta>0$ sufficiently small.

For each $i=1,\ldots,K$, consider the function $g_i(t) := c_i (1 - t^{2^*-2})$. Clearly there exists $\gamma_i>0$ such that $|g_i'(t)|\geq \gamma_i$ for every $t\in [r,R]$. The mean value theorem implies
\[
|g_i(\bar t_i)| = |g_i(\bar t_i) - g_i(1)| \geq \gamma_i |\bar t_i - 1| \qquad \text{for all } \bar t_i \in [r,R].
\]
On the other hand, since $[r,R]^K$ is compact, the map $\bs t\mapsto \prod_{\ell \neq i} t_\ell^2$ is $L$-Lipschitz on $[r,R]^K$ for some $L > 0$. Therefore,
\[
\left|\prod_{\ell \neq i}\bar t_\ell^2 - 1\right| \le L \sum_{\ell=1}^K |\bar t_\ell - 1|.
\]
Combining these two estimates with \eqref{eq:fixed-point-system}, we obtain, for every $i=1,\ldots,K$,
\[
\gamma_i |\bar t_i - 1| \leq d \cdot L \sum_{\ell=1}^K |\bar t_\ell - 1|.
\]
Setting $\gamma := \min_{1 \le i \le K} \gamma_i > 0$ and summing over $i=1,\dots,K$, we get
\[
\gamma \sum_{i=1}^K |\bar t_i - 1| \leq d \cdot K L \sum_{\ell=1}^K |\bar t_\ell - 1|,
\]
or equivalently,
\[
(\gamma - d KL) \sum_{i=1}^K |\bar t_i - 1| \leq 0.
\]
Recall that $d = \lim_{n\to\infty} \beta \int_\Omega \prod_{\ell=1}^K |w_{\ell,n}|^2 dx \leq C \b$ for some $C>0$ independent of $\b$ and of the sequence $\mf u_n$, by definition of $\mf w_n$ and Lemma \ref{lem:uniform-bounds-components}. Therefore, up to reducing $\b$, we deduce that $\bar t_i=1$ for every $i=1,\dots,K$. Since $(t_{1,n}w_{1,n},\ldots,t_{K,n}w_{K,n})\in\mathcal N$, we obtain
\[
\cA
\leq
\lim_{n\to\infty}
I(t_{1,n}w_{1,n},\ldots,t_{K,n}w_{K,n})
=
\lim_{n\to\infty}I(\mf w_n)
=
\frac1N\sum_{i=1}^K b_i.
\]
On the other hand, \eqref{min levels-splitting} yields
\[
\frac1N\sum_{i=1}^K b_i\leq \cB,
\]
thus $\cB\geq \cA$, contradicting \eqref{eq:upper bound B}.
Hence $\mf u\neq\bs0$.

\noindent\textbf{Case 2: $\mf u\neq\bs0$ and at least one component of $\mf u$ vanishes.}

\medskip

After relabelling, assume that
\[
u_1=\cdots=u_m=0,
\qquad
u_i\neq0
\quad\text{for }i=m+1,\ldots,K,
\]
for some $1\leq m<K$. Since at least one component vanishes ($u_1=0$), the interaction term of the limit profile vanishes identically, i.e., $\prod_{i=1}^K |u_i|^2 \equiv 0$. Therefore, for every $i=m+1,\ldots,K$,
\[
-\Delta u_i+\lambda_i u_i
=
\mu_i|u_i|^{2^*-2}u_i,
\]
and hence
\begin{equation}\label{eq:energy-nonzero-components}
    E(\mf u)
    \geq
    \sum_{i=m+1}^K B_{\mu_i}.
\end{equation}

For $i=1,\ldots,m$, we have $b_i>0$ by the same argument as in Case~1. 
We next claim that $b_i=0$ for at least one $i\in\{m+1,\ldots,K\}$. 
Otherwise, $b_i>0$ for every $i=1,\ldots,K$, and we can repeat the 
argument of Case~1. More precisely, we introduce $F_{j,n}$ as before 
and show that $c_j>0$. Indeed, \eqref{min levelsi-ci-d} still holds, 
and if $c_j=0$, then H\"older's inequality implies that $d=0$, 
contradicting the fact that $b_j>0$. The remainder of the argument 
carries over unchanged, yielding
\[
\cA
\leq
\frac1N\sum_{i=1}^K b_i
\leq
E(\mf u)+\frac1N\sum_{i=1}^K b_i
=
\cB,
\]
which again contradicts \eqref{eq:upper bound B}. The claim is proved.

Thus $b_i=0$ for some $i\geq m+1$. Consequently, $w_{i,n}\to0$ strongly in $H_0^1(\Omega)$, and H\"older's inequality yields
\[
\int_\Omega\prod_{\ell=1}^K|w_{\ell,n}|^2\,dx\to0.
\]
Therefore, for every $i=1,\ldots,m$, \eqref{eq:remainder-nehari} gives
\[
\int_\Omega |\nabla w_{i,n}|^2dx
=
\mu_i\int_\Omega |w_{i,n}|^{2^*}dx+o(1).
\]
By the Sobolev inequality,
\[
\int_\Omega |\nabla w_{i,n}|^2dx
\leq
\mu_i S^{-\frac{2^*}{2}}
\left(
\int_\Omega|\nabla w_{i,n}|^2dx
\right)^{\frac{2^*}{2}}
+o(1).
\]
Since $b_i>0$, passing to the limit gives
\[
b_i
\geq
\mu_i^{-\frac{N-2}{2}}S^{\frac N2},
\qquad i=1,\ldots,m.
\]
Combining this with \eqref{min levels-splitting} and \eqref{eq:energy-nonzero-components}, we obtain
\[
\cB
=
E(\mf u)+\frac1N\sum_{i=1}^K b_i
\geq
\sum_{i=m+1}^K B_{\mu_i}
+
\frac1N\sum_{i=1}^m
\mu_i^{-\frac{N-2}{2}}S^{\frac N2},
\]
contradicting \eqref{eq:upper bound B}.

Therefore neither Case~1 nor Case~2 can occur, and hence
\[
u_i\not\equiv0
\qquad\text{for every }i=1,\ldots,K.
\]
Therefore, $\mf u\in\mathcal M$. By the definition of $\cB$,
\[
E(\mf u)\geq \cB,
\]
whereas \eqref{min levels-splitting} gives $E(\mf u)\leq \cB$. Thus
\[
E(\mf u)=\cB.
\]
Therefore, recalling Proposition \ref{prop:minimizers_are_solutions}, we deduce that $\mf u$ is a least energy fully non-trivial solution of \eqref{eq:system 1}.
Finally, replacing each component with its absolute value, we may assume that $\mf u$ is nonnegative and, by the strong maximum principle, strictly positive.
\end{proof}

\section*{AI disclosure}

The authors acknowledge the use of ChatGPT to check the correctness of some proofs and to improve the presentation of the manuscript. All proof strategies and mathematical arguments were developed by the authors. The authors take full responsibility for the entire content of the paper.

\end{document}